\documentclass[fleqn,10pt]{wlscirep}
\usepackage[utf8]{inputenc}
\usepackage[T1]{fontenc}
\usepackage{float}
\usepackage{multirow}
\usepackage{subcaption}
\usepackage{graphicx}
\usepackage{changepage}
\usepackage{tabularx}
\usepackage{lineno}
\title{Fukushima Nuclear Wastewater Discharge: An Evolutionary Game Theory Approach to International and Domestic Interaction and Strategic Decision-Making}

\author[1]{Mingyang Li}
\author[2]{Han Pengsihua}
\author[1,*]{Songqing Zhao}
\author[1,*]{Zejun Wang}
\author[1]{Limin Yang}
\author[1,3]{Weian Liu}
\affil[1]{School of Liberal Arts and Sciences, China University of Petroleum-Beijing at Karamay, Karamay 834000,Xinjiang, China.}
\affil[2]{School of Business Administration, China University of Petroleum-Beijing at Karamay, Karamay 834000,Xinjiang, China.}
\affil[3]{School of Mathematics and Statistics, Wuhan University, Wuhan, Hubei 430072, China.}

\affil[*]{zsq@cup.edu.cn (S.Z.) and wangzj@cupk.edu.cn (Z.W.)}


\begin{abstract}
On August 24, 2023, Japan controversially decided to discharge nuclear wastewater from the Fukushima Daiichi Nuclear Power Plant into the ocean, sparking intense domestic and global debates. This study uses evolutionary game theory to analyze the strategic dynamics between Japan, other countries, and the Japan Fisheries Association. By incorporating economic, legal, international aid, and environmental factors, the research identifies three evolutionarily stable strategies, analyzing them via numerical simulations. The focus is on Japan's shift from wastewater release to its cessation, exploring the myriad factors influencing this transition and their effects on stakeholders' decisions. Key insights highlight the need for international cooperation, rigorous scientific research, public education, and effective wastewater treatment methods. Offering both a fresh theoretical perspective and practical guidance, this study aims to foster global consensus on nuclear wastewater management, crucial for marine conservation and sustainable development.
\end{abstract}
\begin{document}

\flushbottom
\maketitle
%
%
\thispagestyle{empty}

\section*{Introduction}

In 2011, the world witnessed the Fukushima disaster, marking the most severe nuclear catastrophe since the Chernobyl disaster, leaving an indelible imprint on human history. In the initial six weeks following the disaster, the amount of radioactive materials released from Fukushima was detected to be substantial, with the release of cesium amounting to 42\% of the total release in the Chernobyl disaster, and the release of xenon reaching an unprecedented level \cite{1}. In response to this crisis, the Japanese government adopted water injection measures to slow down the melting speed of the fuel rods. However, this measure led to another crisis—large amounts of nuclear-contaminated water needed to be sealed and stored, posing not only a high cost but also a potential permanent burden (International Atomic Energy Agency, IAEA, 2021) \cite{2}. To address this pressing issue, the Japanese government initiated a series of intense discussions and evaluations, exploring various possible solutions. Among numerous solutions, discharging wastewater into the sea was considered a relatively practical approach, albeit filled with controversy and risk. After a series of careful considerations and weighing of pros and cons, on April 13, 2021, the Japanese government officially decided to discharge 1.2 million tons of Fukushima wastewater into the sea \cite{3}, a decision that immediately attracted widespread attention worldwide. Despite the Japanese government's repeated assurances that the discharged polluted water had been treated and the radiation values met safety standards, this action still provoked strong opposition from the international community \cite{4}.

On August 24, 2023, the Japanese government officially commenced the discharge of nuclear wastewater from the Fukushima Daiichi Nuclear Power Station into the ocean, once again triggering extensive international attention and controversy \cite{5}. Although the International Atomic Energy Agency (IAEA) deems Japan's plan to discharge treated water from Fukushima into the sea to be in compliance with international safety standards \cite{6}, many countries have expressed strong opposition due to the uncertainty of the potential dangers and long-term impacts brought by such astonishing discharge. Among them, China firmly opposes Japan's unilateral decision affecting the entire international community and condemns this action \cite{7}. Japan's nuclear wastewater discharge plan is essentially resolving its own nuclear safety threat by sacrificing the safety of the human marine environment \cite{8}. 
As early as the end of 2018, Tokyo Electric Power Company had already admitted that about 80\%, i.e., 890,000 tons, of the then stored 1.1 million tons of treated water still contained \textsuperscript{90}Sr, \textsuperscript{60}Co, \textsuperscript{106}Ru, and many other radioactive nuclides exceeding the limit values, and the purification system failed to adequately remove these radioactive nuclides \cite{9}. 
With the progress of Japan's nuclear wastewater discharge plan, the world is on the verge of experiencing the most severe nuclear pollution. The seas near Fukushima Prefecture are not only the economic source for coastal residents but also an integral part of the Pacific and the world's oceans. The substantial amount of radioactive materials will have immeasurable impacts on marine life, the natural environment, and human health.

The act of discharging has already occurred, and the rational handling of Fukushima's nuclear wastewater is a critical issue that the discharging country and even all stakeholders should focus on. Finding effective approaches is key to resolving governance dilemmas and achieving sustainable development of the oceans. In this context, exploring the long-term relationships between the Japanese government, other countries, and the Japanese Fisheries Association is necessary. Viewing these three parties as stakeholders, the dynamic game among them encompasses the outcomes of handling nuclear contaminated water. We utilize an evolutionary game model to discern the pivotal influencing factors and conduct an in-depth study on the impact of various parameters on the evolutionary process under the condition of Japan refraining from discharging nuclear wastewater, aspiring to contribute to the conservation of the living environment for all of humanity.This paper mainly explores the following three questions:

\begin{itemize}
	\item How to accurately measure the interests and costs of the Japanese government, other countries, and the Japanese Fisheries Association, and how to correctly construct the payoff matrix in the evolutionary game model among the three parties?
	
	\item How to analyze the asymptotic stability of the three parties, and under what conditions can such stability be achieved? What Evolutionarily Stable Strategies (ESS) exist among the three parties?
	
	\item How can an evolutionary game model elucidate the key factors and impacts of parameters in scenarios where Japan refrains from discharging nuclear wastewater, and how can such insights be leveraged for environmental conservation for humanity?
\end{itemize}

In Section 2 of this paper, we review literature pertinent to our study, laying the foundation for the establishment of the theoretical framework and research methodology. Section 3 provides a detailed explanation of each parameter used in the model, offering a clear theoretical basis for subsequent analysis and discussion.
In Section 4, we propose a set of initial assumptions aimed at more accurately constructing the model of benefits and costs for the three parties—Japan, other countries, and the Fisheries Association. Based on these assumptions, we further develop the replicator dynamics equations in the evolutionary game model to reveal the strategic changes and choices of the parties during the game process.
Section 5 delves deeply into the behavioral stability of the three parties under different conditions. We not only analyze various types of asymptotic stability in detail but also explore the necessary conditions and possibilities for achieving such stability, providing important theoretical support for the refinement of the theoretical model and empirical analysis. In Section 6, we validate the feasibility and accuracy of the model through numerical simulations, revealing evolutionary trends under different stable conditions. Specifically, we focus on the scenario where Japan decides to cease ocean discharge, analyzing the impact of key parameters on the evolutionary process under this scenario. This not only provides theoretical references for the three parties to seek optimal strategies in practical problems but also helps in more comprehensively understanding how each party responds to different strategies and environmental variables, thereby affecting the final outcome of the evolutionary game.

The main goal of this paper is to provide scientific theoretical basis and practical guidance for solving the Fukushima nuclear wastewater problem by deeply exploring the dynamic game among the three stakeholders: the Japanese government, other countries, and the Japanese Fisheries Association. We hope that this research can help all parties find the best solution to minimize potential harm to the environment and humanity, achieving sustainable development and harmonious coexistence in human society.

\section*{Literature review}

\subsection*{Chemical hazards of Japan's nuclear wastewater discharge into the sea}
On March 11, 2011, a severe earthquake occurred in the northeastern sea area of Japan, triggering a tsunami and causing serious damage to the Fukushima Daiichi Nuclear Power Plant (FDNPP). The failure of the cooling system led to a temperature rise in the reactor, generating hydrogen and other gases, causing explosions and releasing radioactive gases and fragments into the atmosphere. Simultaneously, untreated cooling water was directly discharged into the northeastern sea area of Japan, making this incident one of the most severe uncontrolled inputs of artificial radionuclides into the ocean in history \cite{cas}. According to the water treatment facility monitoring reports by Tokyo Electric Power Company Holdings (TEPCO, 2012-2020) and the seawater quality monitoring reports by the Nuclear Regulation Authority (NRA, 2013-2020), the radioactive substances mainly include \textsuperscript{3}H, \textsuperscript{14}C, \textsuperscript{134}Cs, \textsuperscript{137}Cs, \textsuperscript{60}Co, \textsuperscript{125}Sb, \textsuperscript{90}Sr, \textsuperscript{129}I, \textsuperscript{99}Tc, \textsuperscript{106}Ru, and \textsuperscript{238}Pu \cite{luy}. Although the concentrations of these isotopes are lower than tritium, they are more likely to integrate into marine organisms and seabed sediments (Buesseler, 2020) \cite{bue}. If nuclear wastewater is discharged into the Pacific Ocean, potential risks will exist for hundreds or even thousands of years in the future. The marine food web tightly connects marine organisms through predation (Albouy et al., 2019) \cite{alb}. Cesium isotopes are the main pollutants in the accumulated seawater inside Japanese nuclear facilities and tanks. Additionally, cesium is a highly soluble radioactive nuclide in seawater, posing long-term radiation risks to the environment and is more likely to integrate into marine organisms or sediments under high concentration factors (Buesseler, 2020) \cite{bue}. The harm of cesium to humans is mainly reflected in its radioactive radiation and human absorption. Moreover, the concentration of \textsuperscript{14}C in treated nuclear wastewater is also high (the Korea Times, 2021b) \cite{the}, it releases low-energy $\beta$ particles and may enter the biosphere and accumulate in marine ecosystems (Williams et al, 2010) \cite{wil}. Co is also a radioactive isotope that releases gamma rays, which can penetrate the human body and cause cell damage (Khajeh et al, 2017) \cite{khaj}. \textsuperscript{90}Sr can simulate calcium in the human body, significantly increasing the risk of osteosarcoma and leukemia (Khani et al., 2012) \cite{khan}. The ocean has a large volume and a complex current system, with strong abilities to dilute and disperse radioactive substances, but long half-life radioactive nuclides will still exist in the marine environment for a long time (Men and Deng et al., 2017) \cite{men}, causing significant consequences to the natural environment, marine life, and human health.

\subsection*{Evolutionary game theory}

Game theory was originally designed to analyze economic behavior and, with modifications, can also be applied to evolving populations. The concept of “Evolutionarily Stable Strategy” (ESS) proposed by Smith (1982) reveals the relevance of the optimal behavior of animals or plants under the influence of the behavior of other individuals \cite{smi1982}. In classical game theory, the basic components of a game include players, strategies, and payoff matrices, and the results of the game can be clearly displayed through the payoff matrices. Evolutionary game theory, originating from classical game theory, has gradually expanded its range of application. It has been widely used in economics, sociology, and statistical physics to analyze various factors affecting the formation of group behavior, producing rich and influential research results (Perc and Jordan et al., 2017; Khoo and Fu et al., 2018; Wang and He et al., 2018)\cite{per, kho, wan, hil}. 

Smith (1973) regarded evolutionary game theory as a theoretical framework for studying the interaction between individual behavior and adaptation. According to this theory, the evolutionary process is a dynamic strategy selection process, where the fitness of an individual depends not only on the chosen strategy but also on interactions with other individuals \cite{smi}. Adami et al. (2016) believed that the main goal in the application of evolutionary game theory is to find suitable strategies or optimal decision sequences to resolve existing conflicts and obtain maximum benefits. However, the competition between different strategies does not occur simultaneously at the same point in time but changes over specific periods. New strategies continuously emerge, compete with existing ones, and may replace the originally dominant strategies by gaining advantages in space and time. Therefore, the success of a strategy depends on its competition with other strategies during a specific period \cite{ada}. Evolutionary game models are important tools for exploring social behavior. Predecessors have constructed various types of evolutionary game models through extensive research to explain the emergence and continuation of various social behaviors. Notably, ESS (Evolutionarily Stable Strategy) is a powerful tool in game theory to find optimal strategies that can prevail in the competition with other strategies (Smith, 1968) \cite{smi1968}.

\subsection*{Application of evolutionary game theory in this study}

Following Japan's decision to discharge nuclear wastewater into the sea, researchers have begun to employ evolutionary game theory to study this issue. 
Su et al. (2023) constructed a Tripartite Evolutionary Game (TEG) model to analyze the evolutionary stability of strategy choices among fishermen, cooperatives, and government departments, exploring influencing factors and their relationships. The results indicate that the profits of cooperatives and the cooperation costs of fishermen are key factors affecting the strategy choices of the three parties in the game. However, this study did not consider the impact of international factors on strategy choices \cite{sum}.
Zhang et al. (2023) adopted a complex network evolutionary game based on the WS small-world model, focusing on the fisherman population, to explore the dynamic evolutionary laws of cooperative behavior diffusion in the development of marine carbon sink fisheries. This study did not involve the roles and strategies of other countries and international organizations, an aspect that is supplemented and expanded upon in this paper \cite{zha}.
Zheng (2022) constructed an evolutionary game model involving fishermen, consumers, and the government, analyzing the impact of each entity's behavior probability on the strategies of other entities and the stability of the entire system. However, this study did not consider the impact of international cooperation and conflicts on strategy evolution \cite{zhe}.
Xin et al. (2022) explored an evolutionary game considering emotional impacts to study the evolutionarily stable strategies of the Japanese government and Japanese fishermen. This study primarily focused on domestic factors and did not fully consider the dynamics and interactions at the international level \cite{xin}.
Xu et al. (2022) established a tripartite evolutionary game model including the International Atomic Energy Agency (IAEA), the Discharging Country (DC), and the Fisheries Cooperative Association of Japan (FCA), aiming to propose management insights for international cooperation. However, this study did not fully explore the changes and evolution of strategies among the parties in the game, an aspect this paper intends to analyze more comprehensively \cite{xul}.

Previous research on Japan's nuclear wastewater discharge was predominantly predictive, given that Japan only commenced the actual discharge on August 24, 2023. Consequently, studies prior to this date were often based on limited information and considerations. In contrast, our research, informed by post-discharge realities, incorporates a broader spectrum of practical factors and employs evolutionary game theory to delve into the strategic decisions and interactions among the stakeholders.
Our approach amalgamates perspectives from economics, environmental science, and international relations, with a particular emphasis on real-time international impacts and responses—elements frequently overlooked in earlier studies. Notably, our work pioneers the exploration of the balance of interests among the Japanese government, other nations, and the Japan Fisheries Association. Furthermore, we are the first to apply evolutionary game theory within this specific context.
Our study not only broadens the knowledge frontier of evolutionary game models in environmental policy and international relations but also, through a tripartite framework, adeptly addresses the literature gap concerning the tripartite balance of interests in Japan's marine discharge policy. In essence, our analysis, rooted in an international dimension, offers a comprehensive and realistic perspective on the dynamic games surrounding Japan's nuclear wastewater issue.

\section*{Symbol definitions}

All symbols and their corresponding meanings are detailed and illustrated in Table~\ref{not}.
\begin{table}[H]
	\centering
	\small 
	\begin{tabular}{@{}l p{12cm}@{}}
		\toprule
		Parameter & Description  \\ \midrule
		\( J \) & Japanese Government\\
		\( C \) & Other Countries \\
		\( F \) & Japanese Fisheries Association \\
		\(x\) & Probability of Japan choosing the discharge strategy \\
		\( 1-x \) & Probability of Japan choosing to cease discharge \\
		\( y \) & Probability of other countries choosing to sanction Japanese discharge \\
		\(1-y\) & Probability of other countries choosing not to sanction Japanese discharge \\
		\( z \) & Probability of the Japanese Fisheries Association opposing the Japanese government \\
		\( 1-z \) & Probability of the Japanese Fisheries Association supporting the Japanese government \\
		\( C_{DJ }\) & Cost of Japanese discharge \\
		\( C_{SJ} \) & Cost of storing nuclear wastewater in Japan \\
		\( I_J \) & International image of Japan \\
		\( C_{IF} \) & International image of the Japanese Fisheries Association \\
		\( C_{LF} \) & Litigation compensation of the Fisheries Association \\
		\( C_{LC} \) & Litigation compensation of other countries \\
		\( E_{RF} \) & Reduction in revenue for the Fisheries Association due to discharge \\
		\( T_{RJ} \) & Reduction in export tax revenue for Japan due to discharge \\
		\( C_{SC} \) & Additional cost for other countries to develop their own seafood products \\
		\( B_{SP}\) & Potential benefits to other countries' seafood industries from introducing Japanese seafood substitutes \\
		\( C_{MJ} \) & Ocean monitoring cost for the Japanese government in the case of discharge \\
		\( C_{MC} \) & Ocean monitoring cost for other countries in the case of discharge \\
		\( C_{HJ} \) & Aid received by Japan from other countries in the case of no discharge \\
		\bottomrule
	\end{tabular}
	\caption{Description of parameters.}
	\label{not}
\end{table}
\section*{Methodology}

\subsection*{Basic assumptions}

\begin{enumerate}[label={\textbf{Assumption} \arabic*:},left=0cm]
	\item In this tripartite game, Japan (J), the National Fisheries Association of Japan (F), and other countries (C) are set as agents. Due to limited information and rationality, the decisions of each party are constrained.
	\item Employing the method of evolutionary game theory, the aforementioned three parties are considered as stakeholders. The dynamic game among these stakeholders involves their respective pursuits of interest and choices of action.
	\item In this game, Japan (J), as the country discharging nuclear-contaminated water, might weigh domestic economic interests against international reputation, deciding whether to continue the discharge of nuclear-contam\\-inated water. The probability of discharging is denoted as \(x\), and conversely, as \(1-x\).
	\item Other countries (C), as additional stakeholders, may exhibit attitudes of either sanctioning or not sanctioning. Some countries might express dissatisfaction with Japan's discharge strategy and adopt sanction measures with a probability of \(y\), to protect international environmental and marine resource interests. Conversely, some might choose not to sanction, with a probability of \(1-y\).
	\item The National Fisheries Association of Japan (F) represents the interests of Japanese fisheries. Their decisions may manifest in two forms: one opposing Japan's discharge strategy with a probability of \(z\), perceiving it as detrimental to fishery resources; the other, with a probability of \(1-z\), might not oppose and thus would endorse Japan's discharge actions.
	\item This study solely analyzes the game strategies of each party from an economic perspective and does not consider the impact of political and other factors on the strategic choices of each participant.
\end{enumerate}

\begin{figure}[H]
		\centering
		\includegraphics[width=15.5cm]{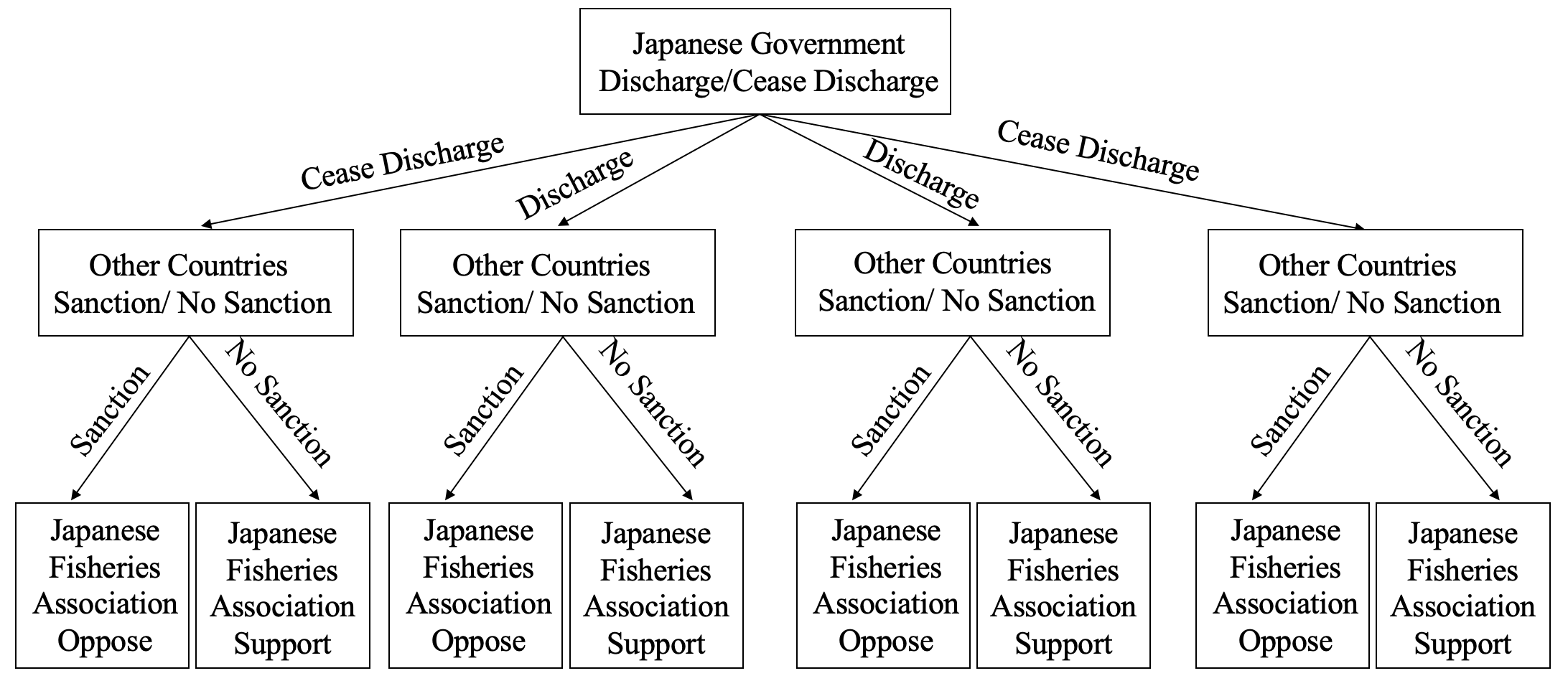}
	\caption{Illustration of all possible scenarios for each participant under different strategies.}
	\label{009}
\end{figure}

Each participant's occurrence in all possible scenarios under different strategies is illustrated in Fig.~\ref{009}.

\subsection*{Evolutionary game model}
The payoff matrix for the Japanese government, other countries, and the Japanese Fisheries Association is shown in Table~\ref{pm}.

\begin{table}[H]
	\centering
	\small
	\caption{Payoff matrix of Japan, other countries, and japanese fisheries association.}
	\resizebox{\textwidth}{!}{ 
		\begin{tabular}{cccccc}  
			\bottomrule
			\multicolumn{4}{c}{\multirow{2}{*}{\textbf{Stakeholders}}} & \multicolumn{2}{c}{\textbf{Japanese Fisheries Association}}  \\
			\cline{5-6}
			\multicolumn{4}{c}{} & Opposition $(z)$ & Acceptance $(1 - z)$ \\
			\hline
			& \multirow{4}{*}{Discharge $(x)$} & & \multirow{2}{*}{Sanction $(y)$} & $
			\begin{aligned}
				&-I_{J}-C_{L F}-C_{L C}-T_{R J}-C_{D J}-C_{M J}, \\
				&-C_{S C}+B_{S P}+C_{L C}-C_{M C}, \\
				&C_{L F}-E_{R F}
			\end{aligned}$ & $
			\begin{aligned}
				&-I_{J}-C_{L C}-T_{R J}-C_{D J}-C_{M J}, \\
				&-C_{S C}+B_{S P}+C_{L C}-C_{M C}, \\
				&-E_{R F}-C_{IF}
			\end{aligned}$ \\  
			& & & & &  \\  
			& & & \multirow{2}{*}{No Sanction $(1-y)$} & $
			\begin{aligned}
				&-C_{L F}-C_{D J}-C_{M J}, \\
				&-C_{M C}, \\
				&C_{L F}
			\end{aligned}$ & $
			\begin{aligned}
				&-C_{D J}-C_{M J}, \\
				&-C_{M C}, \\
				&-C_{IF}
			\end{aligned}$ \\  
			& & & & &  \\  
			& \multirow{4}{*}{No Discharge $(1-x)$ } & & \multirow{2}{*}{Sanction $(y)$} & $
			\begin{aligned}
				&C_{HJ}-C_{SJ}, \\
				&-C_{HJ}, \\
				&0
			\end{aligned}$ & $
			\begin{aligned}
				&C_{HJ}-C_{S J}, \\
				&-C_{HJ}, \\
				&-C_{IF}
			\end{aligned}$\\  
			& & & & & \\ 
			\multirow{-24.5}{*}{\textbf{Japan}} & & \multirow{-24.5}{*}{\textbf{Other Countries}} & \multirow{2}{*}{No Sanction $(1-y)$} & $
			\begin{aligned}
				&-C_{S J}, \\
				&0, \\
				&0
			\end{aligned}$ &  $
			\begin{aligned}
				&-C_{S J}, \\
				&0, \\
				&-C_{IF}
			\end{aligned}$\\  
			& & & &  &  \\  
			\bottomrule
		\end{tabular}	\label{pm}
	} 
\end{table}

Assume that the expected utility of the Japanese government discharging nuclear wastewater is denoted as \(U_{11}\), the expected utility of not discharging is \(U_{12}\), and the average expected utility is \(\bar{U}_{1}\). Then, we have:

\begin{equation}
	\begin{aligned}
		&U_{11}=yz(-I_{J}-C_{LF}-C_{LC}-T_{RJ}-C_{DJ}-C_{MJ})+y(1-z)(-I_{J}-C_{LC}-T_{RJ}-C_{DJ}-C_{MJ})\\
		&+(1-y)z(-C_{LF}-C_{DJ}-C_{MJ})+(1-y)(1-z)(-C_{DJ}-C_{MJ})\\
		&=y(-I_J-C_{LC}-T_{RJ})-zC_{LF}-C_{DJ}-C_{MJ}
	\end{aligned}
\end{equation}

\begin{equation}
	\begin{aligned}
		&U_{12}=yz(C_{HJ}-C_{SJ})+y(1-z)(C_{HJ}-C_{SJ})\\
		&+(1-y)z(-C_{SJ})+(1-y)(1-z)(-C_{SJ})\\
		&=yC_{HJ}-C_{SJ}
	\end{aligned}
\end{equation}

\begin{equation}
	\bar{U}_{1}=xU_{11}+(1-x)U_{12}
\end{equation}

The replicator dynamic equation for the Japanese government, denoted as \(S(x)\), is:

\begin{equation}
	\begin{aligned}
		S(x)=\frac{dx}{dt} & =x(U_{11}-\bar{U}_{1})=x(1-x)(U_{11}-U_{12}) \\
		& =x(1-x)\left[y(-I_J-C_{LC}-T_{RJ}-C_{HJ})-zC_{LF}-C_{DJ}-C_{MJ}+C_{SJ}\right]
	\end{aligned}
	\label{eq1}
\end{equation}

Assume that the expected utility of other countries sanctioning the Japanese government is denoted as \(U_{21}\), the expected utility of other countries not sanctioning the Japanese government is \(U_{22}\), and the average expected utility is \(\bar{U}_{2}\). Then, we have:

\begin{equation}
	\begin{aligned}
		&U_{21}=xz(-C_{SC}+B_{SP}+C_{LC}-C_{MC})+x(1-z)(-C_{SC}+B_{SP}+C_{LC}-C_{MC})\\
		&+(1-x)z(-C_{HJ})+(1-x)(1-z)(-C_{HJ})\\
		&=x(C_{HJ}-C_{SC}+B_{SP}+C_{LC}-C_{MC})-C_{HJ}
	\end{aligned}
\end{equation}

\begin{equation}
	\begin{aligned}
		&U_{22}=xz(-C_{MC})+x(1-z)(-C_{MC})+(1-x)z\cdot0+(1-x)(1-z)\cdot 0\\
		&=-xC_{MC}
	\end{aligned}
\end{equation}

\begin{equation}
	\bar{U}_{2}=yU_{21}+(1-y)U_{22}
\end{equation}

The replicator dynamic equation for other countries, denoted as \(G(y)\), is:

\begin{equation}
	\begin{aligned}
		G(y)=\frac{dy}{dt} & =y(U_{21}-\bar{U}_{2})=y(1-y)(U_{21}-U_{22}) \\
		& =y(1-y)[x(C_{HJ}-C_{SC}+B_{SP}+C_{LC})-C_{HJ}]
	\end{aligned}
	\label{eq2}
\end{equation}

Assume that the expected utility of the Japan Fisheries Association opposing the discharge is denoted as \(U_{31}\), the expected utility of the Japan Fisheries Association supporting the discharge is \(U_{32}\), and the average expected utility is \(\bar{U}_{3}\). Then, we have:

\begin{equation}
	\begin{aligned}
		&U_{31}=xy(C_{LF}-E_{RF})+x(1-y)(C_{LF})+(1-x)y\cdot0+(1-x)(1-y)\cdot0\\
		&=-xyE_{RF}+xC_{LF}
	\end{aligned}
\end{equation}

\begin{equation}
	\begin{aligned}
		&U_{32}=xy(-E_{RF}-C_{IF})+x(1-y)(-C_{IF})+(1-x)y(-C_{IF})+(1-x)(1-y)(-C_{IF})\\
		&=-xyE_{RF}-C_{IF}
	\end{aligned}
\end{equation}

\begin{equation}
	\bar{U}_{3}=zU_{31}+(1-z)U_{32}
\end{equation}

The replicator dynamic equation for the Japan Fisheries Association, denoted as \(P(z)\), is:

\begin{equation}
	\begin{aligned}
		P(z)=\frac{dz}{dt} & =z(U_{31}-\bar{U}_{3})=z(1-z)(U_{31}-U_{32}) \\
		& =z(1-z)(xC_{LF}+C_{IF})
	\end{aligned}
	\label{eq3}
\end{equation}

To further discuss the stable points in evolutionary game theory, we combine the equations~(\ref{eq1})(\ref{eq2})(\ref{eq3}) to obtain the system of replicator dynamic equations as follows:
\begin{equation}
	\left\{
	\begin{aligned}
		S(x) &= \frac{d x}{d t}  = x(1-x)\left[y(-I_j-C_{LC}-T_{RJ}-C_{HJ})-zC_{LF}-C_{DJ}-C_{MJ}+C_{SJ}\right]\\
		G(y) &=\frac{d y}{d t} = y(1-y)[x(C_{HJ}-C_{SC}+B_{SP}+C_{LC})-C_{HJ}] \\
		P(z) &= \frac{d z}{d t}  = z(1-z)(xC_{LF}+C_{IF})\\
	\end{aligned}
	\right.
	\label{fuzhi}
\end{equation}

\section*{Evolutionary stability analysis}
\subsection*{Asymptotic stability analysis of the three parties}

Let the dynamic replicator equations~(\ref{fuzhi}) be defined as \(S(x)=G(y)=P(z)=0\). By solving these, we can identify eight pure strategy equilibria and one mixed strategy equilibrium, which are \(\gamma_{1}(0,0,0)^{T}\), \(\gamma_{2}(1,0,0)^{T}\), \(\gamma_{3}(0,1,0)^{T}\), \(\gamma_{4}(0,0,1)^{T}\), \(\gamma_{5}(1,1,0)^{T}\), \(\gamma_{6}(1,0,1)^{T}\), \(\gamma_{7}(0,1,1)^{T}\), \(\gamma_{8}(1,1,1)^{T}\), and \(\gamma_{9}\left(x^{*}, y^{*}, z^{*}\right)^{T}\), where \(\gamma_{9}\left(x^{*}, y^{*}, z^{*}\right)^{T}\) is determined by the following equations:

\begin{equation}
	\left\{
	\begin{aligned}
		& y^*(-I_j-C_{LC}-T_{RJ}-C_{HJ})-z^*C_{LF}-C_{DJ}-C_{MJ}+C_{SJ}=0\\
		& x^*(C_{HJ}-C_{SC}+B_{SP}+C_{LC})-C_{HJ}=0\\
		&x^*C_{LF}+C_{IF}=0\\
	\end{aligned}
	\right.
\end{equation}

In asymmetric games, Evolutionarily Stable Strategies (ESS) are pure strategy equilibria, hence we only consider the asymptotic stability of pure strategies. The stability of equilibrium points in evolutionary games is typically determined by the stability of the Jacobian matrix, as illustrated below.

\begin{equation}
	J=\left[\begin{array}{lll}
		J_{11} & J_{12} & J_{13} \\
		J_{21} & J_{22} & J_{23} \\
		J_{31} & J_{32} & J_{33}
	\end{array}\right]=\left[\begin{array}{lll}
		\frac{\partial S(x)}{\partial x} & \frac{\partial S(x)}{\partial y}  & \frac{\partial S(x)}{\partial z}  \\
		\frac{\partial G(y)}{\partial x} &\frac{\partial G(y)}{\partial y}   & \frac{\partial G(y)}{\partial z}\\
		\frac{\partial P(z)}{\partial x}& \frac{\partial P(z)}{\partial y} & \frac{\partial P(z)}{\partial z}
	\end{array}\right]
\end{equation}

where
\begin{equation}
	J_{11}=(1-2x)\left[y(-I_j-C_{LC}-T_{RJ}-C_{HJ})-zC_{LF}-C_{DJ}-C_{MJ}+C_{SJ}\right] 
\end{equation}
\begin{equation}
	J_{12}=x(1-x)(-I_j-C_{LC}-T_{RJ}-C_{HJ}) 
\end{equation}
\begin{equation}
	J_{13}=x(1-x)(-C_{LF})
\end{equation}
\begin{equation}
	J_{21}=y(1-y)(C_{HJ}-C_{SC}+B_{SP}+C_{LC})
\end{equation}
\begin{equation}
	J_{22}=(1-2y)[x(C_{HJ}-C_{SC}+B_{SP}+C_{LC})-C_{HJ}] 
\end{equation}
\begin{equation}
	J_{23}= 0 
\end{equation}
\begin{equation}
	J_{31}= z(1-z)C_{LF}
\end{equation}
\begin{equation}
	J_{32}=0 
\end{equation}
\begin{equation}
	J_{33}=(1-2z)(xC_{LF}+C_{IF})
\end{equation}

According to Lyapunov stability theory, the stability of equilibrium points in evolutionary games can be determined by the eigenvalues of the Jacobian matrix\cite{theo}. An equilibrium point is asymptotically stable, and the strategy is an Evolutionarily Stable Strategy (ESS), only if all the eigenvalues of the Jacobian matrix are negative.

By analyzing the computed results for each point (in Appendix), the stability analysis of these points can be found as shown in Table~\ref{st}.

\begin{table}[H]
	\centering
	\footnotesize
	\caption{Stability analysis of equilibrium point.}
	\label{st}
	\newcolumntype{C}{>{\centering\arraybackslash}X}
		\begin{tabular}{ccccc}
	\hline 
	\multirow{2}{*}{Equilibrium Point} & \multicolumn{2}{c}{Eigenvalues} & \multirow{2}{*}{Stability} & \multirow{2}{*}{Condition} \\
	\cline{2-3}
	& $\lambda(\lambda_{1}, \lambda_{2}, \lambda_{3})$ & Symbol & & \\
	\hline
	$\gamma_{1}(0,0,0)^{T}$ & 
	\begin{tabular}{@{}l@{}}
		$C_{SJ} - C_{MJ} - C_{DJ}$ \\
		$-C_{HJ}$ \\
		$C_{IF}$
	\end{tabular} &  (*,-, +) & Non-ESS & \\
	
	$\gamma_{2}(1,0,0)^{T}$ & 
	\begin{tabular}{@{}l@{}}
		$C_{DJ} + C_{MJ} - C_{SJ}$ \\
		$B_{SP} + C_{LC} - C_{SC}$ \\
		$C_{LF} + C_{IF}$
	\end{tabular} &  (*,*, +) & Non-ESS& \\
	
	$\gamma_{3}(0,1,0)^{T}$ & 
	\begin{tabular}{@{}l@{}}
		$C_{SJ} - C_{HJ} - C_{LC} - C_{MJ} - C_{DJ} - I_{J} - T_{RJ}$ \\
		$C_{HJ}$ \\
		$C_{IF}$
	\end{tabular} &  (*,+, +) & Non-ESS & \\
	
	$\gamma_{4}(0,0,1)^{T}$ & 
	\begin{tabular}{@{}l@{}}
		$C_{SJ} - C_{LF} - C_{MJ} - C_{DJ}$ \\
		$-C_{HJ}$ \\
		$-C_{IF}$
	\end{tabular} &  (-,-, -) & ESS & \begin{tabular}{@{}l@{}}
		$C_{{SJ}}<\left({C}_{{LF}}+{C}_{{MJ}}+{C}_{{DJ}}\right)$
	\end{tabular} \\	
	
	$\gamma_{5}(1,1,0)^{T}$ & 
	\begin{tabular}{@{}l@{}}
		$C_{DJ} + C_{HJ} + C_{LC} + C_{MJ} - C_{SJ} + I_{J} + T_{RJ}$ \\
		$C_{SC} - C_{LC} - B_{SP}$ \\
		$C_{LF} + C_{IF}$
	\end{tabular} &  (*,*, +) & Non-ESS & \\
	
	$\gamma_{6}(1,0,1)^{T}$ & 
	\begin{tabular}{@{}l@{}}
		$C_{DJ} + C_{LF} + C_{MJ} - C_{SJ}$ \\
		$B_{SP} + C_{LC} - C_{SC}$ \\
		$-C_{IF} - C_{LF}$
	\end{tabular} &  (-,-,-) & ESS & 		\begin{tabular}{@{}l@{}}
		$C_{SJ} > (C_{LF} + C_{MJ} + C_{DJ})$ \\
		$B_{SP} + C_{LC} < C_{SC}$ \\
	\end{tabular}
	\\
	
	$\gamma_{7}(0,1,1)^{T}$ & 
	\begin{tabular}{@{}l@{}}
		$C_{SJ} - C_{HJ} - C_{LC} - C_{LF} - C_{MJ} - C_{DJ} - I_{J} - T_{RJ}$ \\
		$C_{HJ}$ \\
		$-C_{IF}$
	\end{tabular} &  (*,+,-) &  Non-ESS& \\
	
	$\gamma_{8}(1,1,1)^{T}$ & 
	\begin{tabular}{@{}l@{}}
		$C_{DJ} + C_{HJ} + C_{LC} + C_{LF} + C_{MJ} - C_{SJ} + I_{J} + T_{RJ}$ \\
		$C_{SC} - C_{LC} - B_{SP}$ \\
		$-C_{IF} - C_{LF}$
	\end{tabular} &  (-,-, -) &  ESS &		\begin{tabular}{@{}l@{}}
		$C_{{SJ}}>{C}_{{DJ}}+{C}_{{HJ}}+{C}_{{LC}}+{C}_{{LF}}+{C}_{{MJ}}+{I}_{{J}}+{T}_{{RJ}} $\\
		${C_{SC}}<{C_{LC}}+{B_{SP}}$
	\end{tabular} \\
		\bottomrule
	\end{tabular}
	\begin{flushleft}
		\footnotesize{* Note: "-" indicates that the eigenvalue of the Jacobian matrix is negative, "+" means that the eigenvalue is positive, and "*" indicates that the sign of the eigenvalue is uncertain.}
	\end{flushleft}
\end{table}

Table~\ref{st} shows that the equilibrium points $\gamma_{1}(0,0,0)^{T}$ , 	$\gamma_{2}(1,0,0)^{T}$ , $\gamma_{3}(0,1,0)^{T}$ , $\gamma_{5}(1,1,0)^{T}$ and $\gamma_{7}(0,1,1)^{T}$  all have positive eigenvalues, hence they are unstable points. The equilibrium points $\gamma_{4}(0,0,1)^{T}$ , $\gamma_{6}(1,0,1)^{T}$  and $\gamma_{8}(1,1,1)^{T}$  can all have negative eigenvalues. Therefore, they are stable points.
\begin{itemize}
	\item When $\mathrm{C}_{\mathrm{SJ}} < \left(\mathrm{C}_{\mathrm{LF}} + \mathrm{C}_{\mathrm{MJ}} + \mathrm{C}_{\mathrm{DJ}}\right)$, the eigenvalues of the Jacobian matrix corresponding to the evolutionary stable point $\gamma_{4}(0,0,1)^{T}$ are all negative, indicating that $\gamma_{4}(0,0,1)^{T}$ is an evolutionarily stable point. Under this condition, when the cost for Japan to store nuclear wastewater is less than the sum of the costs of discharging it into the sea, marine monitoring fees, and litigation compensation to the Fisheries Association, Japan is more likely to choose the "non-discharge" strategy. By adopting this strategy, Japan can reduce the direct impact on the marine ecosystem, respond to public concerns, protect the image and reputation of the government, potentially gain a competitive advantage in the field of environmental protection technology, and promote international cooperation. Other countries tend to choose the "non-sanction" strategy, as they might worry that sanctioning Japan when it chooses "non-discharge" could have negative impacts on bilateral trade, investment, and economic cooperation, thereby damaging their own economic interests, leading them to choose a more conservative strategy. For the Japanese Fisheries Association, it tends to choose the "oppose" strategy because, even under this condition where the cost of storing nuclear wastewater is less than other related costs, storing nuclear wastewater cannot completely eliminate the potential risks to the marine environment. The association needs to consider that such risks might cause long-term damage to the marine ecosystem and fishery resources, bringing uncertainty to fishing activities and the livelihoods of fishermen.
	
	\item When $C_{SJ} > (C_{LF} + C_{MJ} + C_{DJ})$ and $B_{SP} + C_{LC} < C_{SC}$, the eigenvalues of the Jacobian matrix corresponding to the evolutionary stable point $\gamma_{6}(1,0,1)^{T}$ are all negative, indicating that $\gamma_{6}(1,0,1)^{T}$ is an evolutionarily stable point. At this time, the cost for Japan to store nuclear wastewater exceeds the sum of the costs of discharging it into the sea, marine monitoring fees, and litigation compensation to the Japanese Fisheries Association; the potential benefits and litigation compensation obtained by other countries from introducing substitutes for Japanese seafood are less than the additional costs for other countries to develop their own seafood products. In this scenario, Japan is more likely to choose the "discharge" strategy, allowing it to reduce the costs associated with storing nuclear wastewater. Other countries are more inclined to choose the "non-sanction Japan" strategy, as they can fill the domestic supply gap of Japanese seafood by introducing substitutes and increase their own seafood exports, hence they might prefer introducing seafood substitutes over pursuing litigation compensation. The Japanese Fisheries Association may choose to "oppose" due to a sharp decline in seafood income caused by marine pollution.
	
	\item When $\mathrm{C}_{\mathrm{SJ}} > \left(\mathrm{C}_{\mathrm{DJ}}+\mathrm{C}_{\mathrm{HJ}}+\mathrm{C}_{\mathrm{LC}}+\mathrm{C}_{\mathrm{LF}}+\mathrm{C}_{\mathrm{MJ}}+\mathrm{I}_{\mathrm{J}}+\mathrm{T}_{\mathrm{RJ}}\right)$ and $\mathrm{C_{SC}} < \mathrm{C_{LC}}+\mathrm{B_{SP}}$, all the eigenvalues of the Jacobian matrix corresponding to the evolutionary stable point $\gamma_{8}(1,1,1)^{T}$ are negative. This denotes that $\gamma_{8}(1,1,1)^{T}$ is an evolutionarily stable point. Under such circumstances, the cost incurred by Japan to store nuclear wastewater surpasses the aggregate of the costs related to sea discharge, acquiring aid from other nations, compensations from litigations to other nations and the domestic Fisheries Association, marine surveillance fees, impacts on Japan's global image, and the decrement in export tax revenue due to sea discharge. Concurrently, the supplementary cost for other nations to cultivate their seafood is inferior to the aggregate of the litigation compensations and the prospective gains to their seafood industry from adopting Japanese seafood substitutes. In this context, Japan is inclined to opt for the "discharge" strategy, mitigating the expenses linked to nuclear wastewater storage. The Japanese Fisheries Association is predisposed to "oppose," as its cardinal objective is to safeguard and advance the domestic fisheries sector's interests. The discharge of nuclear wastewater into the sea could inflict irreversible damages to the marine ecosystem and pose potential hazards to fisheries resources and the fishing sector. This discharge may contaminate fisheries resources and compromise the quality of the catches, adversely impacting the revenue and living conditions of the fishermen.  When the additional cost for other nations to cultivate their seafood is less than the combined litigation compensations and the prospective advantages to their seafood industry from adopting Japanese seafood substitutes, other nations are likely to select the "sanction" strategy. This choice aims to preserve the interests of their seafood sector, secure economic gains from litigation compensations, and sustain a commendable international image and reputation.
\end{itemize}

\section*{Numerical simulation}
\subsection*{Numerical simulation of different stable points}
In the aforementioned tripartite evolutionary game model, the strategy equilibrium of each game participant is influenced by the strategy choices of other participants. To more intuitively explore the strategy selection process of the three stable equilibrium points mentioned above, we set different parameter values according to different scenarios for numerical simulation to analyze the evolutionary trajectories of stable strategies of the Japanese government, other countries, and the Japanese Fisheries Association.

As can be seen from Table~\ref{st}, the stable points are $\gamma_{4}(0,0,1)^{T}$, $\gamma_{6}(1,0,1)^{T}$, and $\gamma_{8}(1,1,1)^{T}$. To delve deeper into the dynamic characteristics and evolutionary processes of these stable points, we simulate the evolutionary trajectories of these three stable points. We have set a group of parameter values that meet the stable conditions for each stable point. These parameter values will be used to simulate the evolutionary trajectories of these stable points to understand their dynamic behavior more accurately. The specific parameter values are shown in Table~\ref{shuzhi}, with the initial points all set to $[0.5,0.5,0.5]$.

\begin{table}[H]
	\centering
	\caption{Basic Parameter Settings for Three Stable Equilibrium Points.}
	\small
	\label{shuzhi}
	\begin{tabular}{ccccccccccccc}
		\toprule
		Condition&Stable point & $I_{J}$ & $C_{LC}$ & $T_{RJ}$ & $C_{HJ}$ & $C_{LF}$ & $C_{DJ}$ & $C_{MJ}$ & $C_{SJ}$ & $C_{IF}$ &    $B_{SP}$ &  $C_{SC}$  \\
		\midrule
		Condition 1&$\gamma_{4}(0,0,1)^{T}$  & 20 & 8 & 5 & 10 & 35 & 3 & 6 & 30 & 1 & 1 & 30 \\
		
		Condition 2&$\gamma_{6}(1,0,1)^{T}$ & 20 & 8 & 5 & 10 & 20 & 3 & 6 & 30 & 1 & 1 & 30\\
		
		Condition 3&$\gamma_{8}(1,1,1)^{T}$ & 20 & 10 & 5 & 10 & 5 & 3 & 6 & 80 & 1 & 1 & 10  \\
		\bottomrule
	\end{tabular}
\end{table}

\subsubsection*{Evolutionary trajectory under Condition 1}
\begin{figure}[H]
	\centering
	\begin{subfigure}{0.45\linewidth}
		\includegraphics[width=\linewidth]{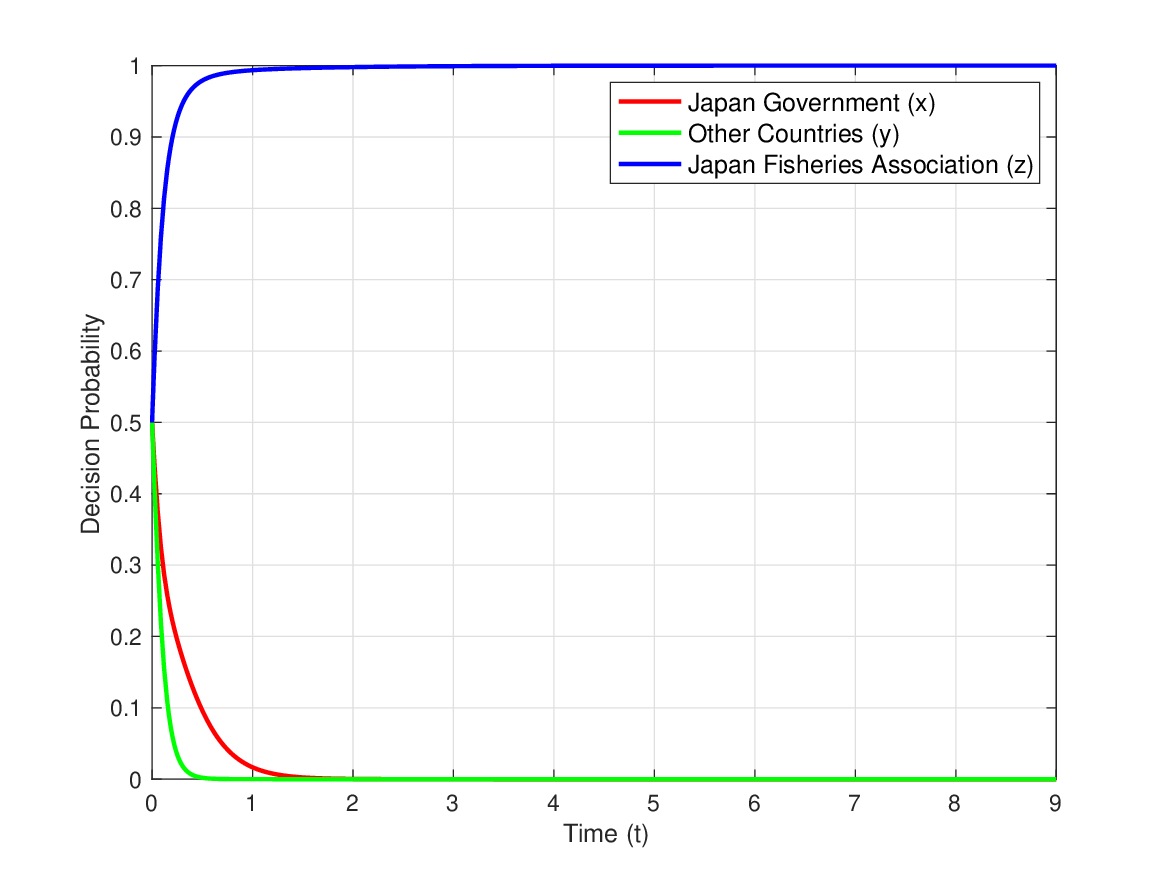}
		\caption{Evolution pattern over time.}
		\label{1a}
	\end{subfigure}
	\begin{subfigure}{0.45\linewidth}
		\includegraphics[width=\linewidth]{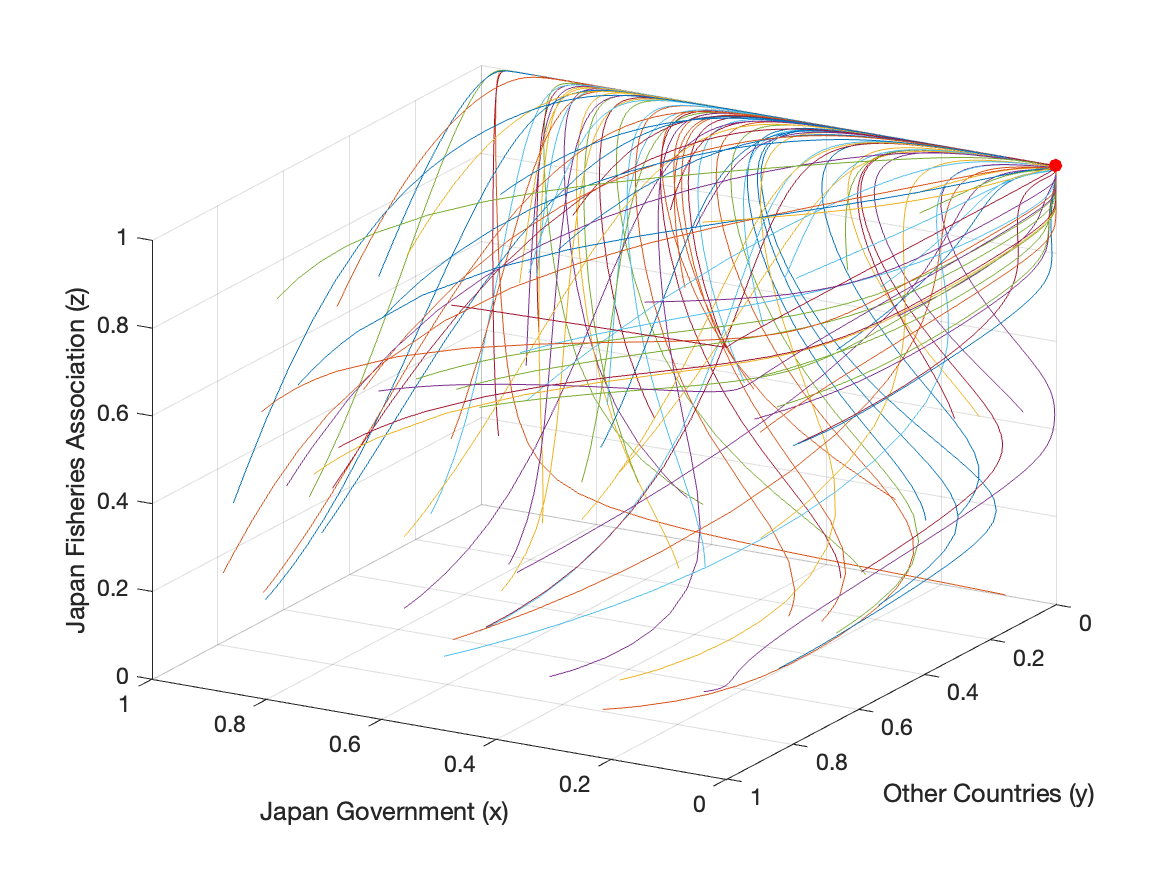}
		\caption{Evolutionary trajectory.}
		\label{1b}
	\end{subfigure}
	\caption{Evolutionary trajectory concerning $\gamma_{4}(0,0,1)^{T}$ under Condition 1.}
	\label{1}
\end{figure}

When the system satisfies Condition 1, i.e., $\mathrm{C}_{\mathrm{SJ}}<\left(\mathrm{C}_{\mathrm{LF}}+\mathrm{C}_{\mathrm{MJ}}+\mathrm{C}_{\mathrm{DJ}}\right)$, we refer to the parameters of Condition 1 in Table~\ref{shuzhi} for analysis. This condition implies that when the cost for Japan to store nuclear wastewater is less than the total of other related costs (these costs include the cost of discharging into the sea, marine monitoring fees, and litigation compensation to the Japanese domestic Fisheries Association), the Japanese government is more likely to choose the "non-discharge" strategy. This is because, under such circumstances, choosing to store nuclear wastewater instead of discharging it into the sea is more economical for Japan from a financial perspective.
The results of numerical simulation further confirm this viewpoint. The simulation results show that under this condition, the stable point of the system is $(0,0,1)$. Fig.~\ref{1a} displays the evolution of the probability of each decision-maker over time, clearly showing how each entity gradually tends towards this stable point over time. Fig.~\ref{1b} illustrates the overall evolutionary trajectory under Condition 1, providing us with a macroscopic perspective to observe the dynamic behavior of the system.

\subsubsection*{Evolutionary trajectory under Condition 2}
\begin{figure}[H]
	\centering
	\begin{subfigure}{0.45\linewidth}
		\includegraphics[width=\linewidth]{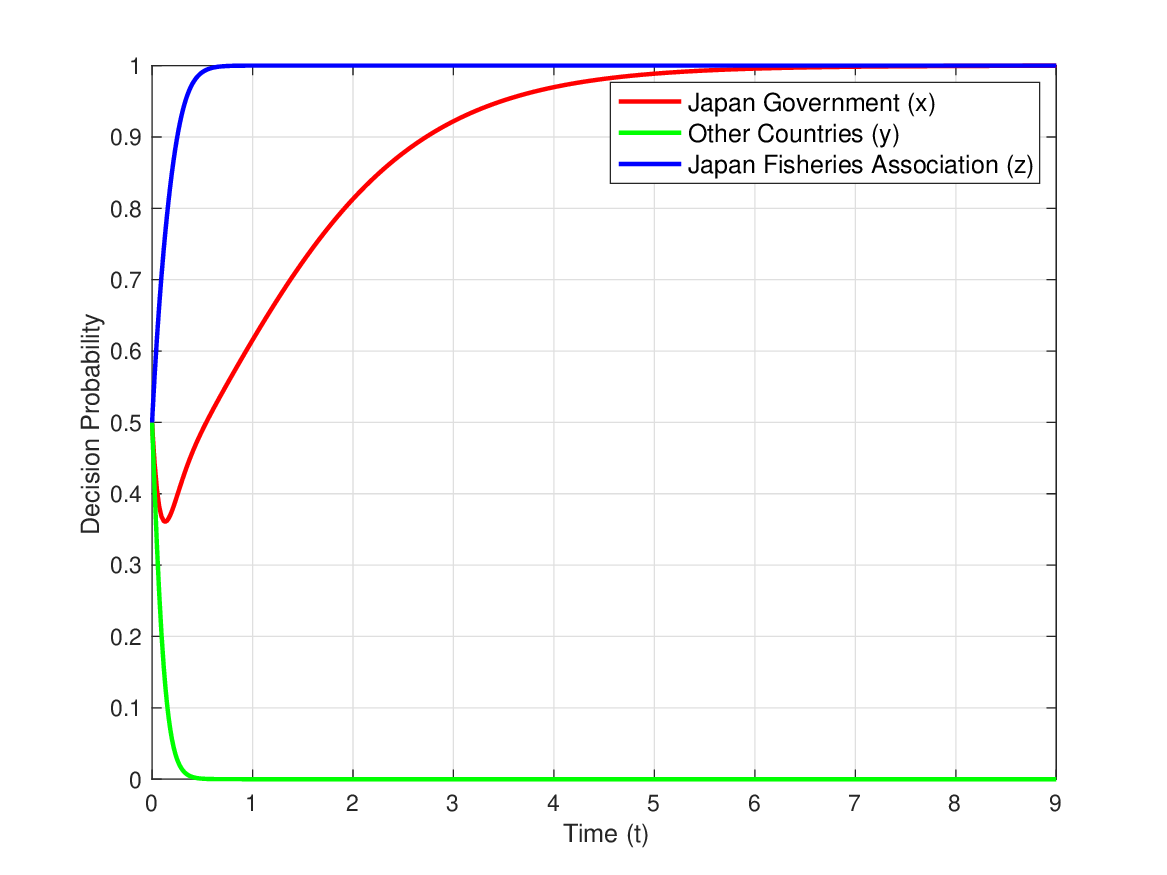}
		\caption{Evolution of probability over time.}
		\label{2a}
	\end{subfigure}
	\begin{subfigure}{0.45\linewidth}
		\includegraphics[width=\linewidth]{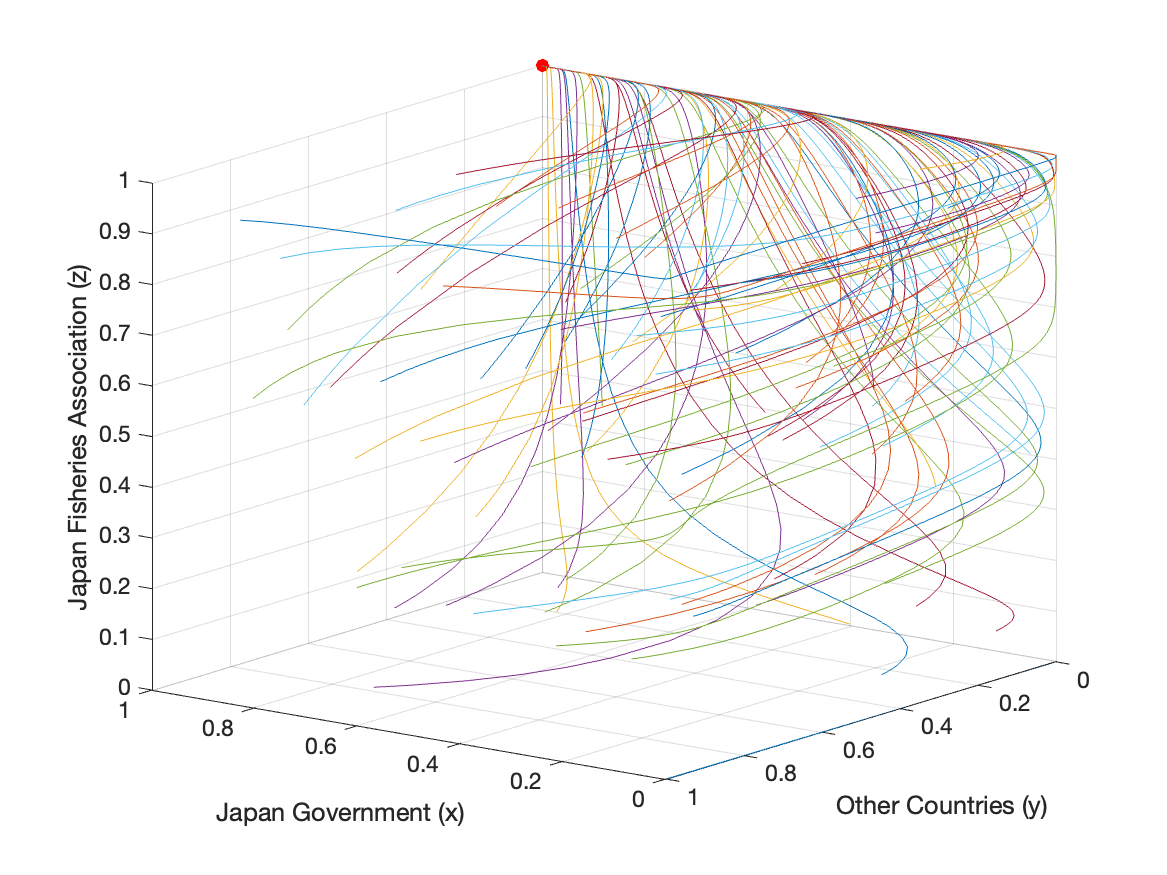}
		\caption{Evolutionary trajectory.}
		\label{2b}
	\end{subfigure}
	\caption{Evolutionary trajectory concerning $\gamma_{6}(1,0,1)^{T}$ under Condition 2.}
	\label{fig:overall_label}
\end{figure}

When the system satisfies Condition 2, i.e., $C_{SJ} > (C_{LF} + C_{MJ} + C_{DJ})$ and $B_{SP} + C_{LC} < C_{SC}$, we refer to the parameters of Condition 3 in Table~\ref{shuzhi} for analysis. This condition implies that the cost for Japan to store nuclear wastewater has exceeded the costs of discharging into the sea, marine monitoring fees, and litigation compensation to the Japanese domestic Fisheries Association; the potential benefits and litigation compensation obtained by other countries through the introduction of Japanese seafood substitutes are less than the additional costs of developing their own seafood. Under these circumstances, Japan is more likely to choose the "discharge" strategy. This decision involves considerations from multiple aspects, including economic costs, environmental impacts, social responsibility, and public opinion.
The results of numerical simulation further confirm this viewpoint. The simulation results show that under this condition, the stable point of the system is $(1,0,1)$.Fig.~\ref{2a} displays the evolution of the probability of each decision-maker over time, clearly showing how each entity gradually tends towards this stable point over time. Fig.~\ref{2b} illustrates the overall evolutionary trajectory under Condition 2, providing us with a macroscopic perspective to observe the dynamic behavior of the system.

\subsubsection*{Evolutionary trajectory under Condition 3}

\begin{figure}[H]
	\centering
	\begin{subfigure}{0.45\linewidth}
		\includegraphics[width=\linewidth]{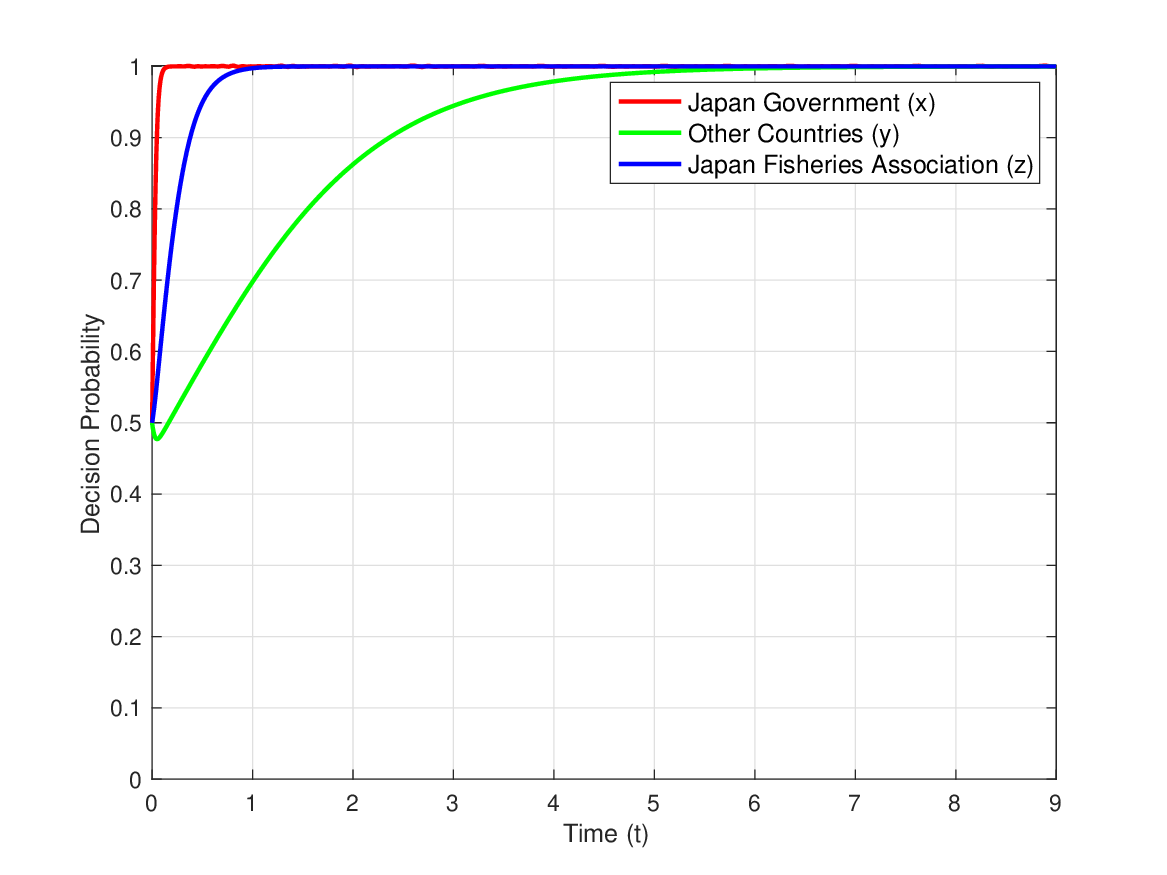}
		\caption{Evolution of probability over time.}
		\label{3a}
	\end{subfigure}
	\begin{subfigure}{0.45\linewidth}
		\includegraphics[width=\linewidth]{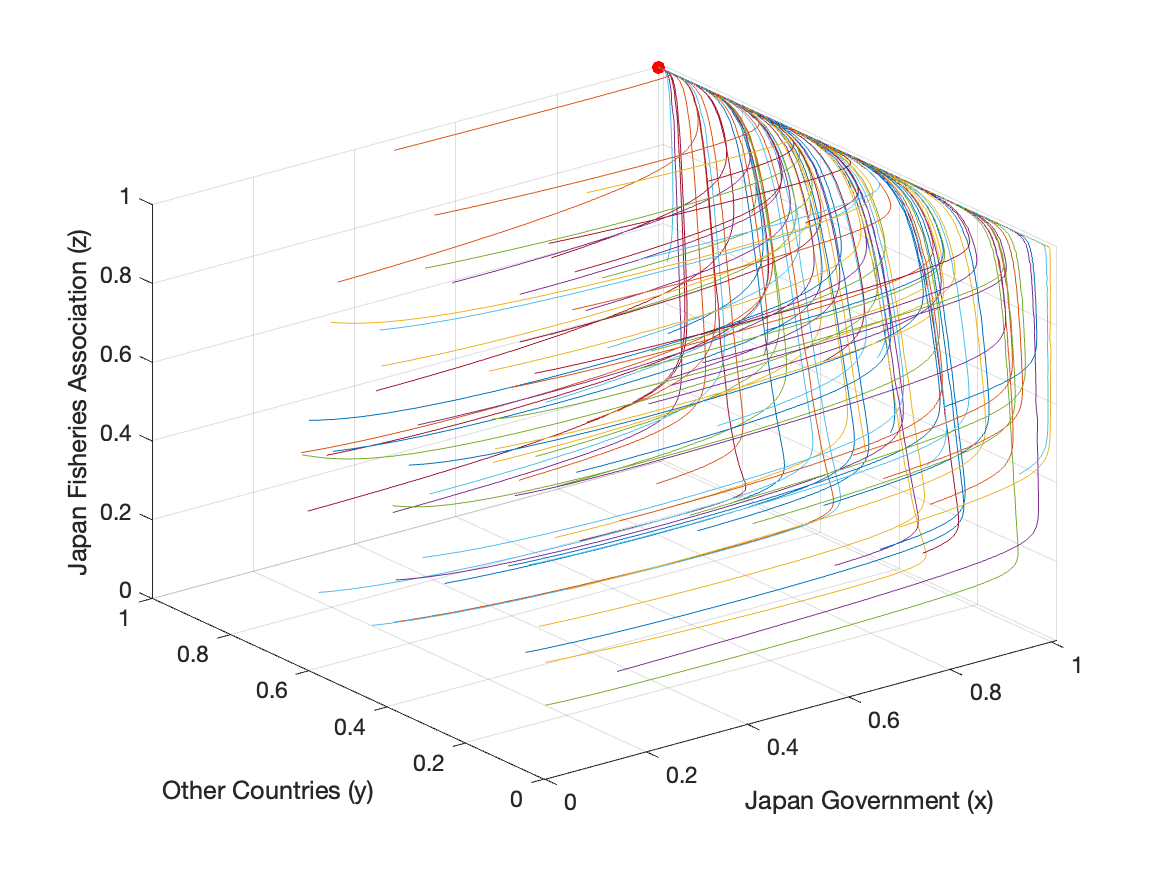}
		\caption{Evolutionary trajectory.}
		\label{3b}
	\end{subfigure}
	\caption{Evolutionary trajectory concerning $\gamma_{8}(1,1,1)^{T}$ under Condition 3.}
	\label{fig:overall_label}
\end{figure}

When the system satisfies Condition 3, i.e., $\mathrm{C}_{\mathrm{SJ}} > \left(\mathrm{C}_{\mathrm{DJ}}+\mathrm{C}_{\mathrm{HJ}}+\mathrm{C}_{\mathrm{LC}}+\mathrm{C}_{\mathrm{LF}}+\mathrm{C}_{\mathrm{MJ}}+\mathrm{I}_{\mathrm{J}}+\mathrm{T}_{\mathrm{RJ}}\right)$ and $\mathrm{C_{SC}} < \mathrm{C_{LC}}+\mathrm{B_{SP}}$, we refer to the parameters of Condition 3 in Table~\ref{shuzhi} for analysis. The cost for Japan to store nuclear wastewater has already exceeded the costs of discharging into the sea, marine monitoring fees, loss of international image, reduction in export taxes, and litigation compensation to the Japanese domestic Fisheries Association and other countries; the additional cost for other countries to develop their own seafood is less than the potential benefits and litigation compensation obtained by introducing Japanese seafood substitutes. Under these circumstances, Japan is more likely to choose the "discharge" strategy.
The results of numerical simulation further confirm this viewpoint. The simulation results show that under this condition, the stable point of the system is $(1,0,1)$. Fig.~\ref{3a} displays the evolution of the probability of each decision-maker over time, clearly showing how each entity gradually tends towards this stable point over time. Fig.~\ref{3b} illustrates the overall evolutionary trajectory under Condition 3, providing us with a macroscopic perspective to observe the dynamic behavior of the system.

\subsection*{Impact of key parameters on evolutionary trajectories}
Our goal is for Japan to ultimately choose the strategy to stop discharging into the sea, i.e., Condition 1. To delve deeper into the key parameters under this condition, such as the cost of Japan discharging into the sea (\(C_{DJ}\)), the cost of Japan storing nuclear wastewater (\(C_{SJ}\)), litigation compensation from other countries (\(C_{LC}\)), litigation compensation from the Fisheries Association (\(C_{LF}\)), aid received by Japan from other countries (\(C_{HJ}\)), and the initial strategy choice probabilities of Japan, other countries, and the Japanese Fisheries Association (\(x_0\), \(y_0\), \(z_0\)) on the Evolutionarily Stable Strategy (ESS) of the three parties, this paper has conducted a series of numerical simulations under Condition 1.

\begin{table}[H]
	\centering
	\caption{Initial Parameter Settings under Condition 1.}
	\label{tab:initial-parameters}
	\begin{tabular}{lccccccccccc}
		\toprule
		Stable point & $I_{J}$ & $C_{LC}$ & $T_{RJ}$ & $C_{HJ}$ & $C_{LF}$ & $C_{DJ}$ & $C_{MJ}$ & $C_{SJ}$ & $C_{IF}$ & $B_{SP}$ & $C_{SC}$ \\
		\midrule
		$\gamma_{4}(0,0,1)^{T}$ & 20 & 8 & 5 & 10 & 30 & 3 & 6 & 30 & 1 & 1 & 30 \\
		\bottomrule
	\end{tabular}
\end{table}

Under these settings, we explore the impact of these parameters on the evolution results and trajectories of the discharge strategy of the Japanese government, the sanction strategy of other countries, and the opposition strategy of the Japanese Fisheries Association, with specific parameter values as shown in Table~\ref{tab:initial-parameters}. Our simulation results reveal how these key parameters influence the strategy choices and evolutionary dynamics of each party in the game, providing in-depth insights for our understanding of the formulation and implementation of international environmental policies.

\subsubsection*{Cost of Japan's discharge into the sea $C_{DJ}$}
\begin{figure}[H]
	\centering
	\begin{subfigure}[t]{0.32\linewidth}
		\includegraphics[width=\linewidth]{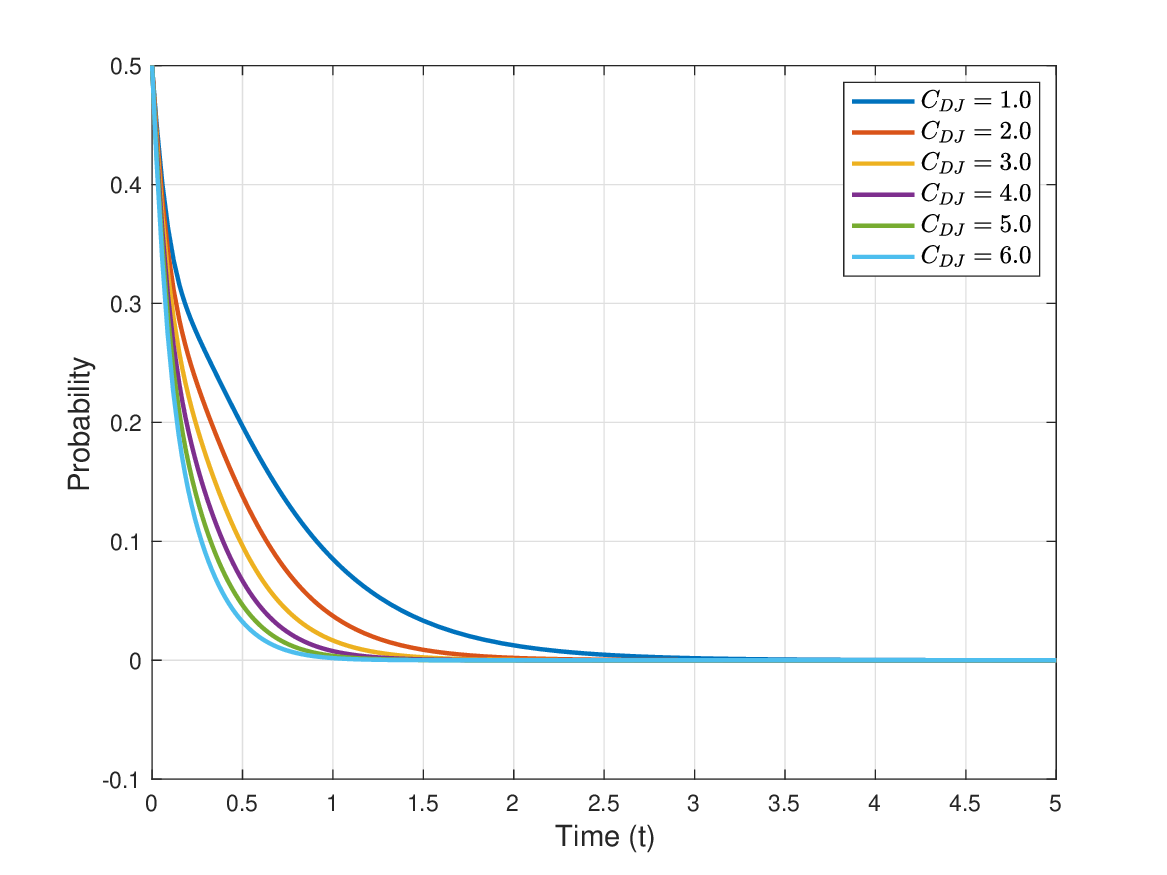}
		\caption{Evolution of the Japanese Government over Time.}
		\label{f1a}
	\end{subfigure}
	\begin{subfigure}[t]{0.32\linewidth}
		\includegraphics[width=\linewidth]{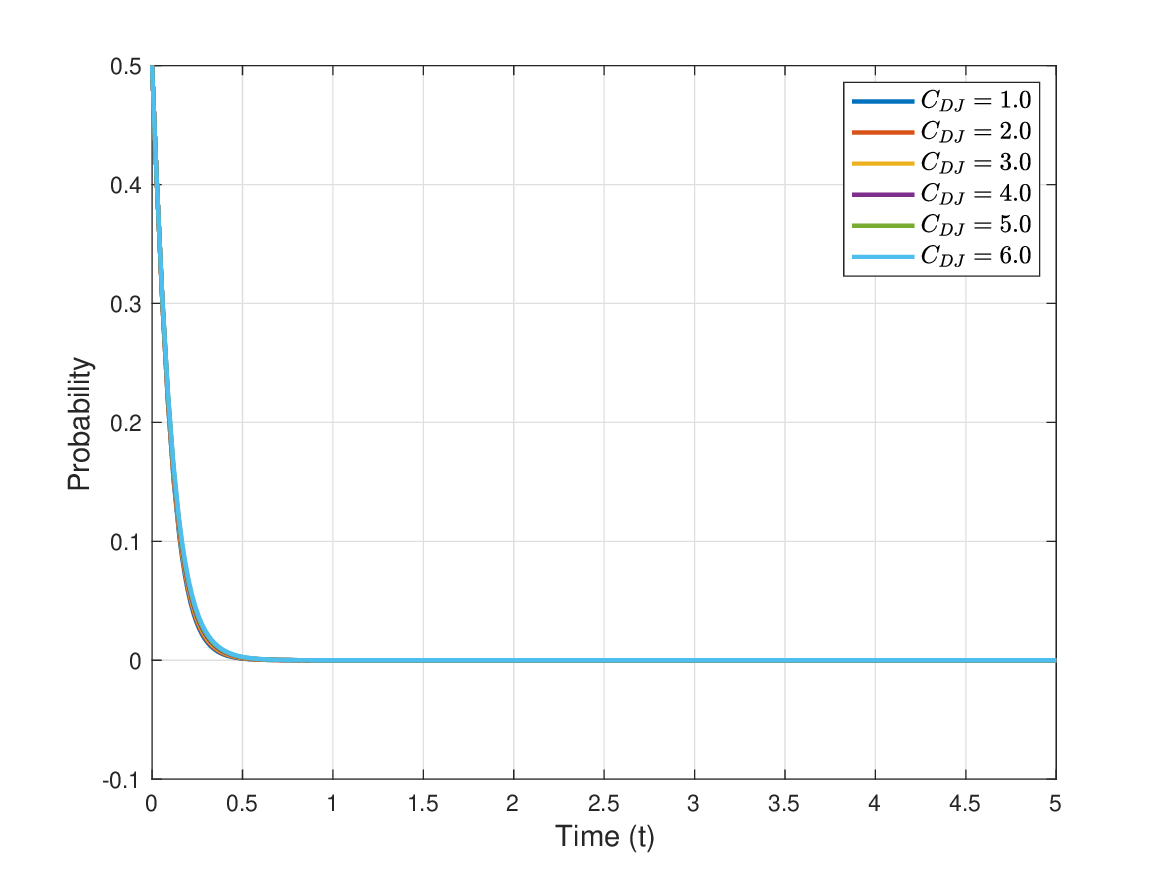}
		\caption{Evolution of Other Countries over Time.}
		\label{f1b}
	\end{subfigure}
	\begin{subfigure}[t]{0.32\linewidth}
		\includegraphics[width=\linewidth]{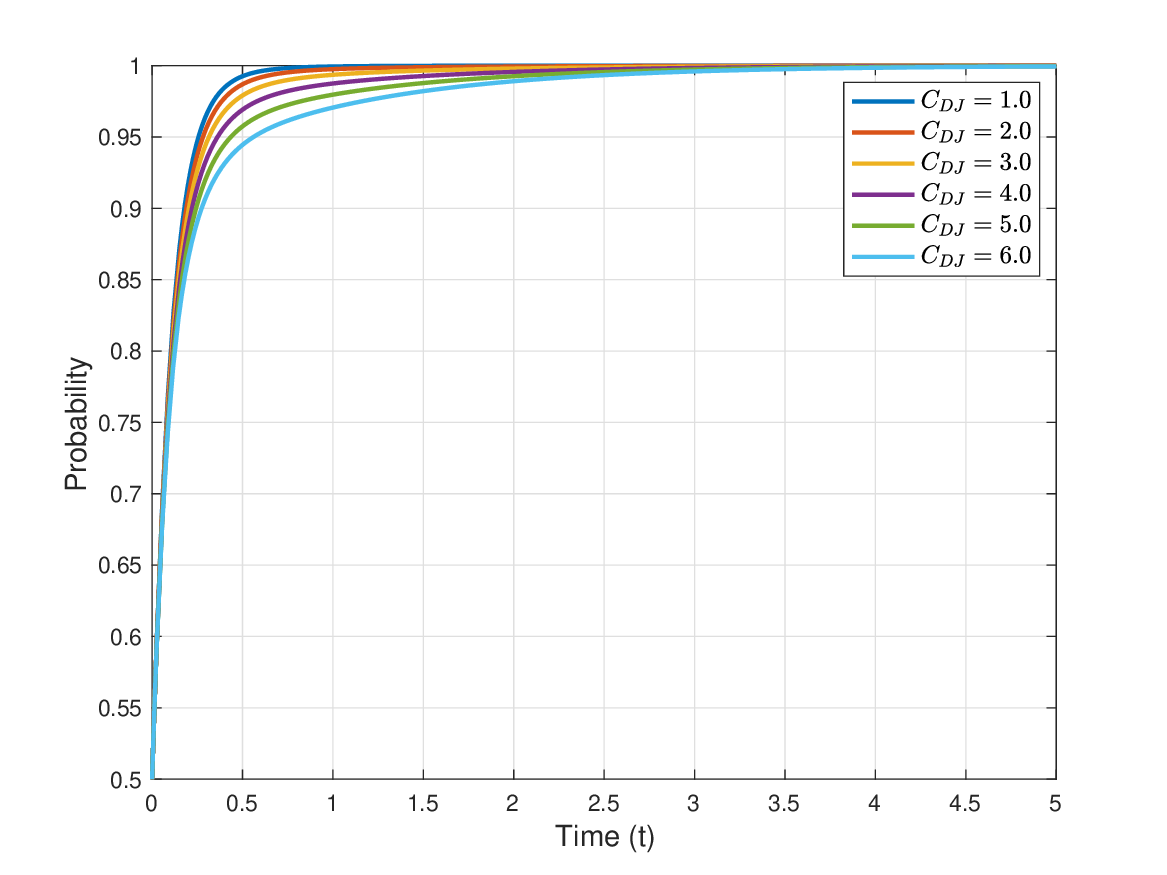}
		\caption{Evolution of the Japanese Fisheries Association over Time.}
		\label{f1c}
	\end{subfigure}
	\caption{Impact of Japan's Discharge Cost $C_{DJ}$ on Evolution}
	\label{f1}
\end{figure}

We explore the cost of Japan's discharge into the sea, denoted as\(C_{DJ}\), on the evolutionary outcomes and trajectories of the three game entities. We assume \(C_{DJ} = 1.0, 2.0, 3.0, 4.0, 5.0\), and \(6.0\),  while maintaining other parameters as shown in Table~\ref{tab:initial-parameters}. The specific impacts of varying \(C_{DJ}\) are illustrated in Fig.~\ref{f1}.
As illustrated in Fig.~\ref{f1a}, with the increase in \(C_{DJ}\), the evolutionary stable point of the Japanese government, \(\gamma_4(0, 0, 1)^T\), remains constant, but the evolutionary speed rapidly accelerates, meaning the Japanese government reaches the equilibrium point more swiftly. This might imply that higher discharge costs prompt the Japanese government to quickly solidify its strategy to mitigate additional costs induced by uncertainties. 
As depicted in Fig.~\ref{f1b}, with the increase in \(C_{DJ}\), the evolutionary stable point of other countries, \(\gamma_4(0, 0, 1)^T\), also remains constant, and the evolutionary speed does not significantly change. This might be attributed to the decision-making of other countries being predominantly influenced by their domestic factors and international political elements, having little association with Japan’s discharge costs. 
As shown in Fig.~\ref{f1c}, with the increase in \(C_{DJ}\), the evolutionary stable point of the Japanese Fisheries Association, \(\gamma_4(0, 0, 1)^T\), likewise remains constant, but the evolutionary speed decelerates, meaning the Japanese Fisheries Association reaches the equilibrium point more slowly. This might be because the Japanese Fisheries Association requires more time to assess and adapt to higher discharge costs and the potential impacts of these costs on fisheries and related industries.

In conclusion, higher costs of processing nuclear wastewater not only augment the probability of Japan adopting the discharge strategy but also expedite its evolutionary speed. The evolutionary speed and strategy selection of other countries and the Japanese Fisheries Association, on the other hand, are more influenced by their respective internal and external factors.

\subsubsection*{Cost of storing nuclear wastewater in Japan \(C_{SJ}\)}
\begin{figure}[H]
	\centering
	\begin{subfigure}[t]{0.32\linewidth}
		\includegraphics[width=\linewidth]{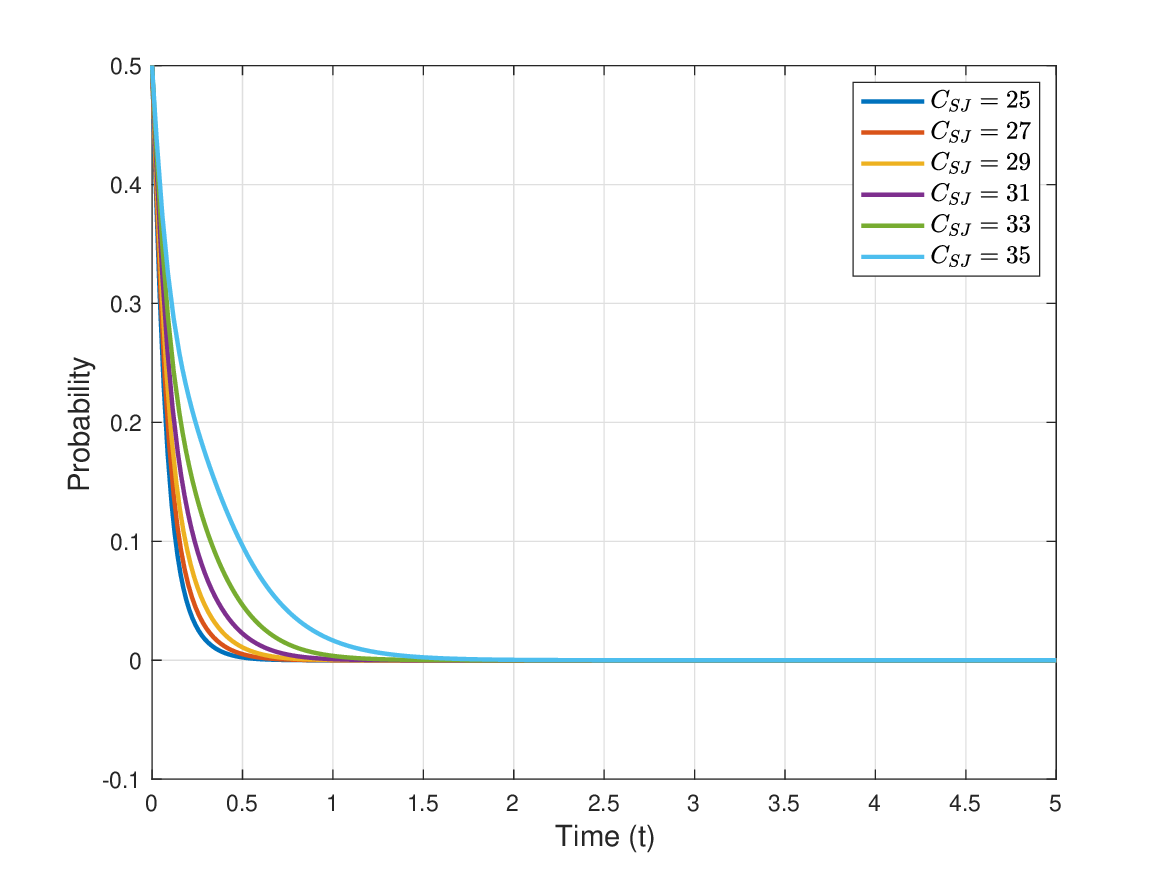}
		\caption{Evolution of the Japanese Government over Time.}
		\label{f2a}
	\end{subfigure}
	\begin{subfigure}[t]{0.32\linewidth}
		\includegraphics[width=\linewidth]{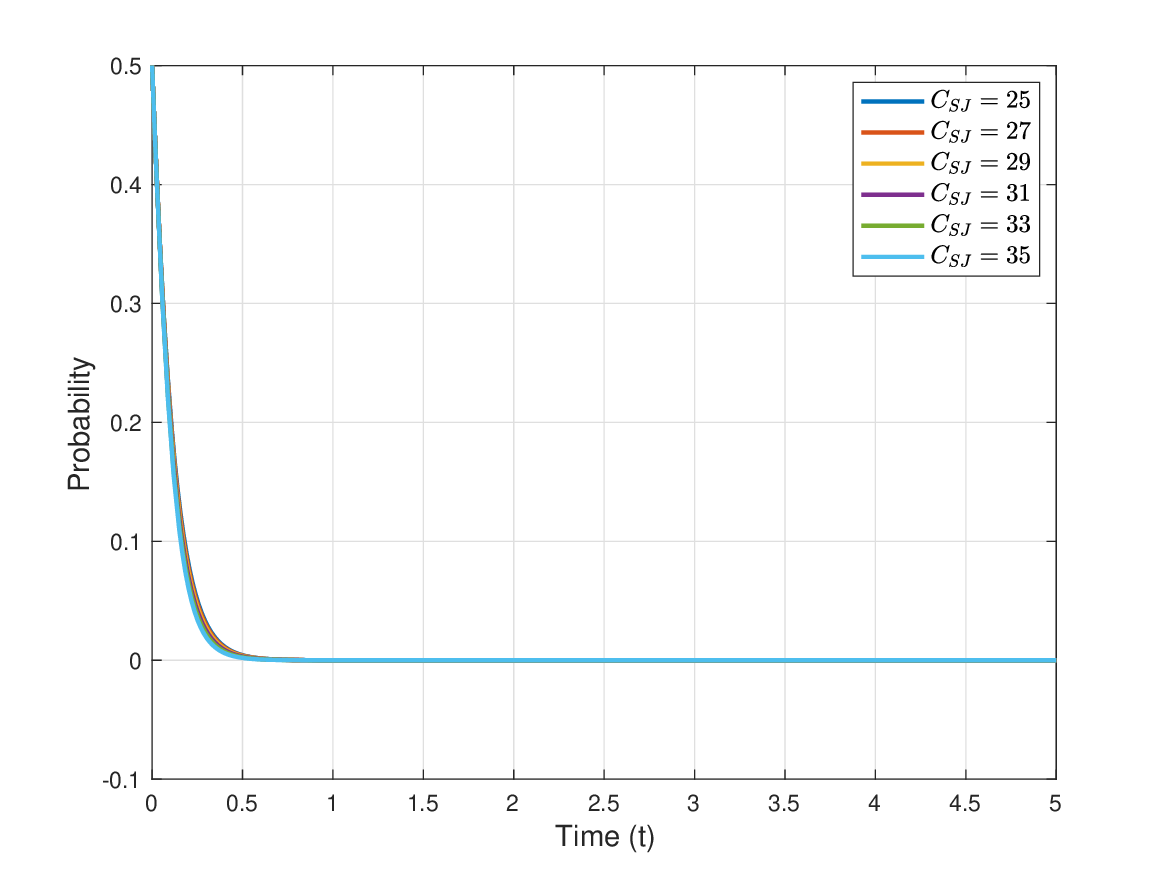}
		\caption{Evolution of Other Countries over Time.}
		\label{f2b}
	\end{subfigure}
	\begin{subfigure}[t]{0.32\linewidth}
		\includegraphics[width=\linewidth]{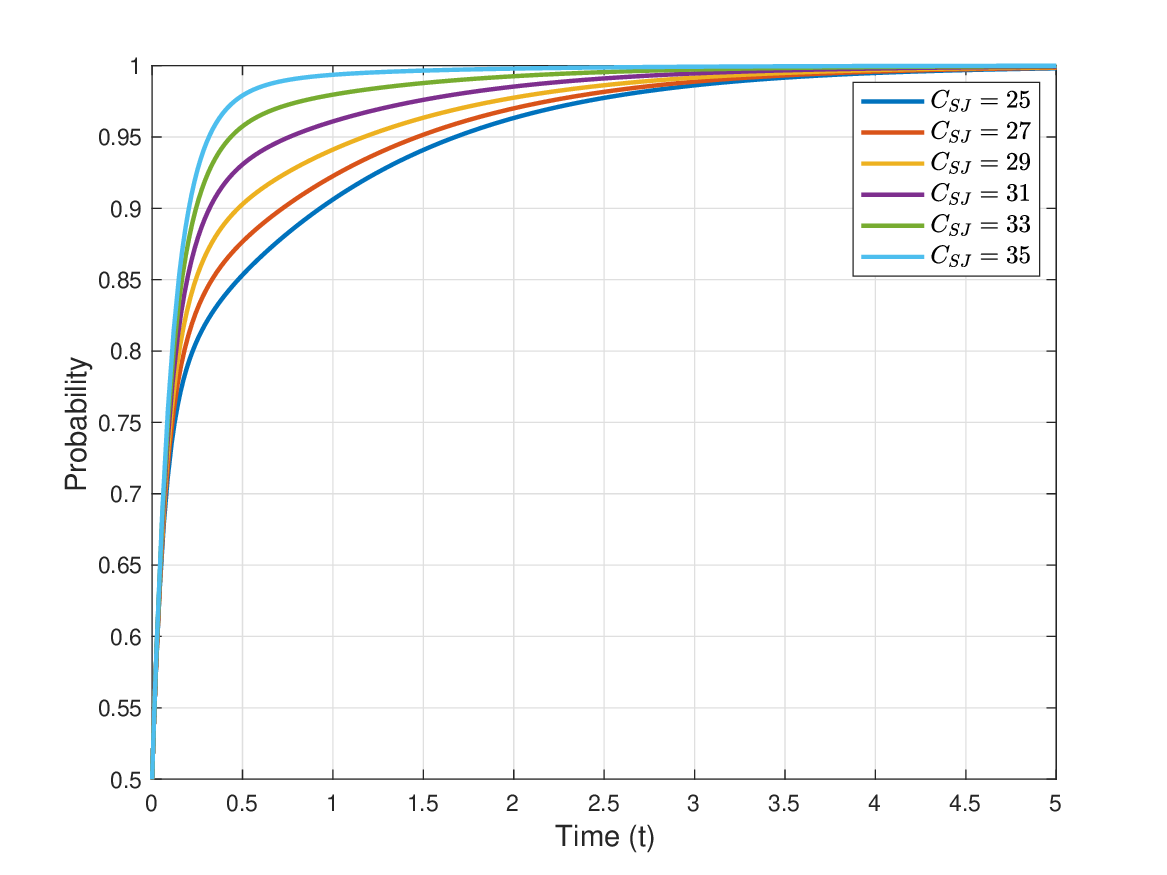}
		\caption{Evolution of the Japanese Fisheries Association over Time.}
		\label{f2c}
	\end{subfigure}
	\caption{Impact of the Cost of Storing Nuclear Wastewater in Japan $C_{SJ}$ on Evolution.}
	\label{f2}
\end{figure}

We explore the cost of storing nuclear wastewater in Japan, denoted as\(C_{SJ}\), on the evolutionary outcomes and trajectories of the three game entities. We assume \(C_{SJ} = 25, 27, 29, 31, 33,\) and \(35\), while maintaining other parameters as shown in Table~\ref{tab:initial-parameters}. The specific impacts of varying \(C_{SJ}\) are illustrated in Fig.~\ref{f2}.
As \(C_{SJ}\) increases, the evolutionary stable point of the Japanese government, \(\gamma_4(0, 0, 1)^T\), remains constant, but the evolutionary speed decreases, implying that the Japanese government requires more time to reach the equilibrium point. This may suggest that higher costs of storing nuclear wastewater necessitate a longer deliberation period for the Japanese government to formulate its strategy, aiming to mitigate the additional costs brought about by uncertainties. As depicted in Fig.~\ref{f2b}, with the increase in \(C_{SJ}\), the evolutionary stable point of other countries, \(\gamma_4(0, 0, 1)^T\), also remains unchanged, and the evolutionary speed is essentially constant. This could be attributed to the fact that the decisions of other countries are predominantly influenced by their domestic considerations and international political factors, having little correlation with Japan's cost of discharging into the sea. As shown in Fig.~\ref{f2c}, with the increase in \(C_{SJ}\), the evolutionary stable point of the Japanese Fisheries Association, \(\gamma_4(0, 0, 1)^T\), also remains constant, but the evolutionary speed consistently accelerates, meaning the Japanese Fisheries Association reaches the equilibrium point more rapidly. This indicates that the Japanese Fisheries Association places a higher emphasis on environmental issues and desires the Japanese government to adopt stricter environmental protection measures promptly.

In conclusion, the varying costs of storing nuclear wastewater in Japan not only influence the strategic deliberation and evolutionary pace of the Japanese government but also reflect the environmental priorities and expectations of the Japanese Fisheries Association, while other countries remain relatively unaffected, maintaining their strategic stances based on their internal and international considerations.

\subsubsection*{Litigation compensation to other countries \(C_{LC}\)}

\begin{figure}[H]
	\centering
	\begin{subfigure}[t]{0.32\linewidth}
		\includegraphics[width=\linewidth]{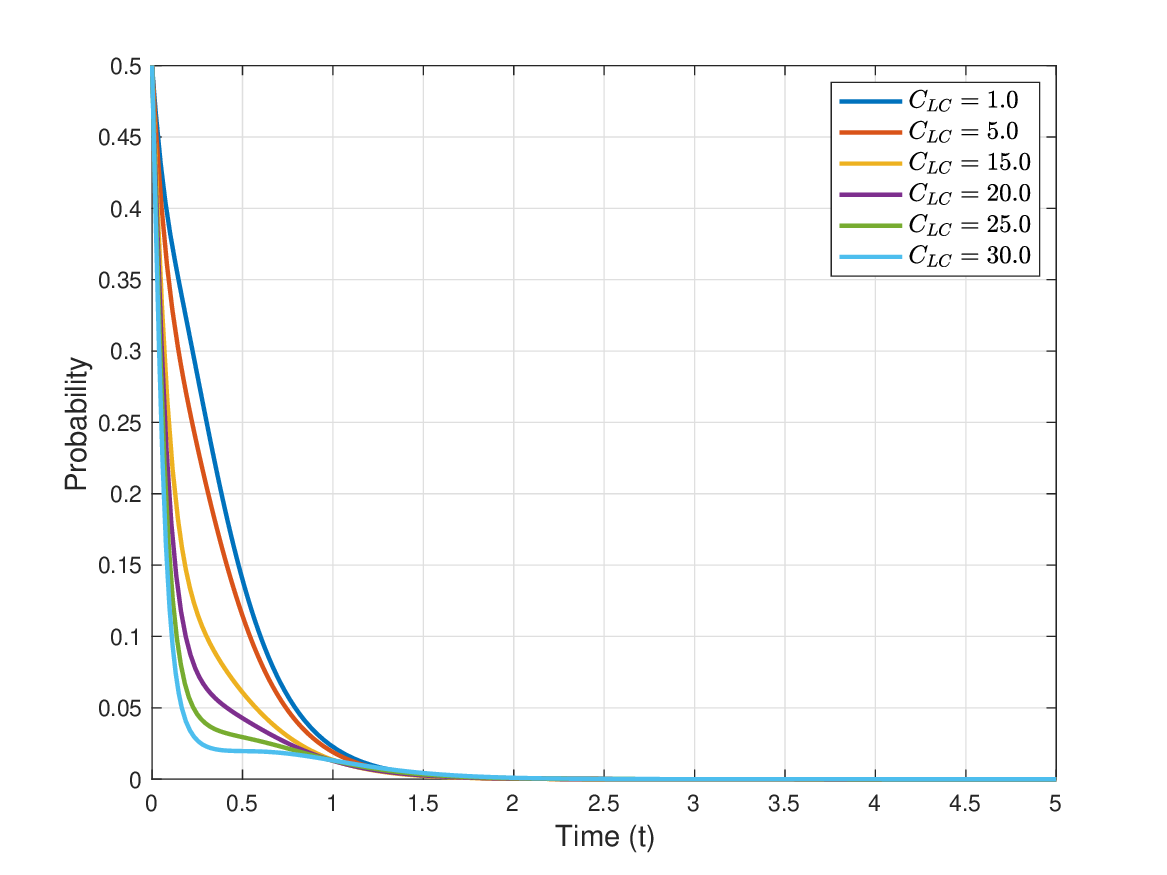}
		\caption{Evolution of the Japanese Government over Time.}
		\label{f3a}
	\end{subfigure}
	\begin{subfigure}[t]{0.32\linewidth}
		\includegraphics[width=\linewidth]{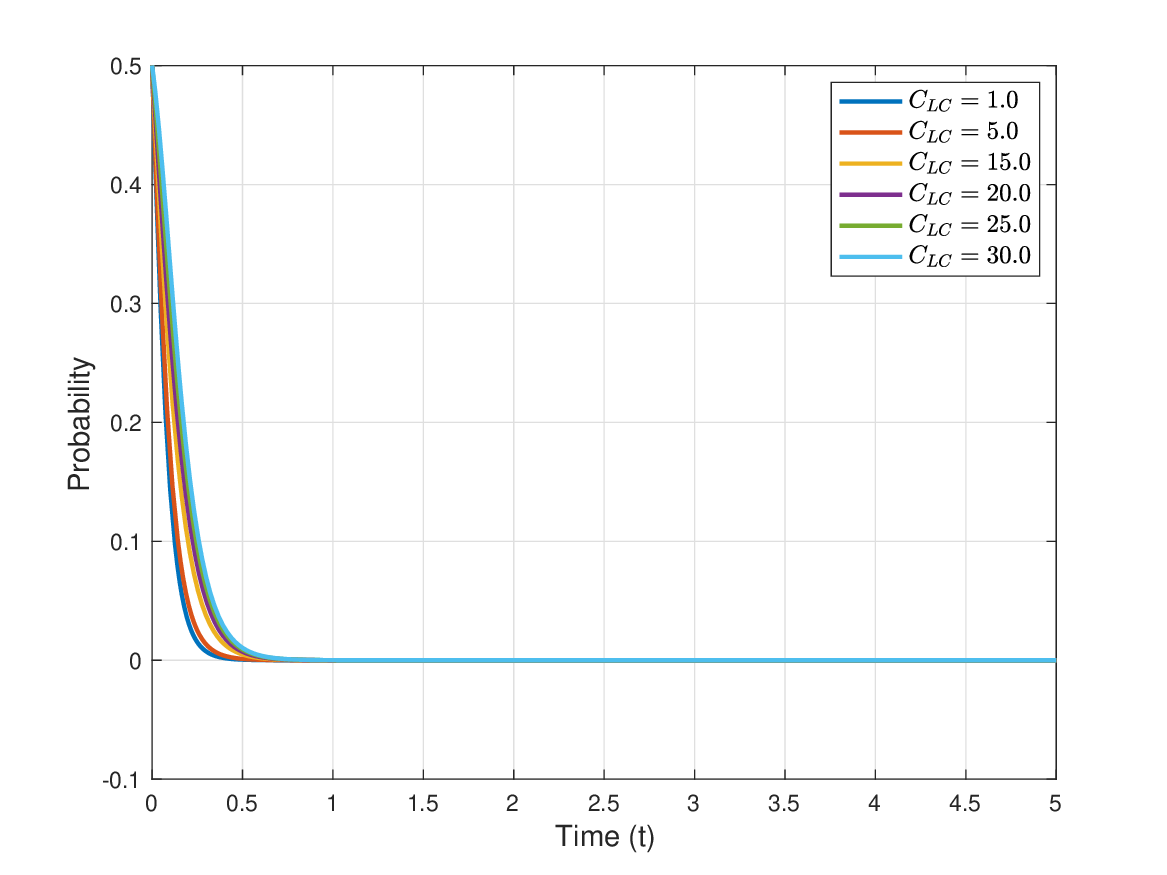}
		\caption{Evolution of Other Countries over Time.}
		\label{f3b}
	\end{subfigure}
	\begin{subfigure}[t]{0.32\linewidth}
		\includegraphics[width=\linewidth]{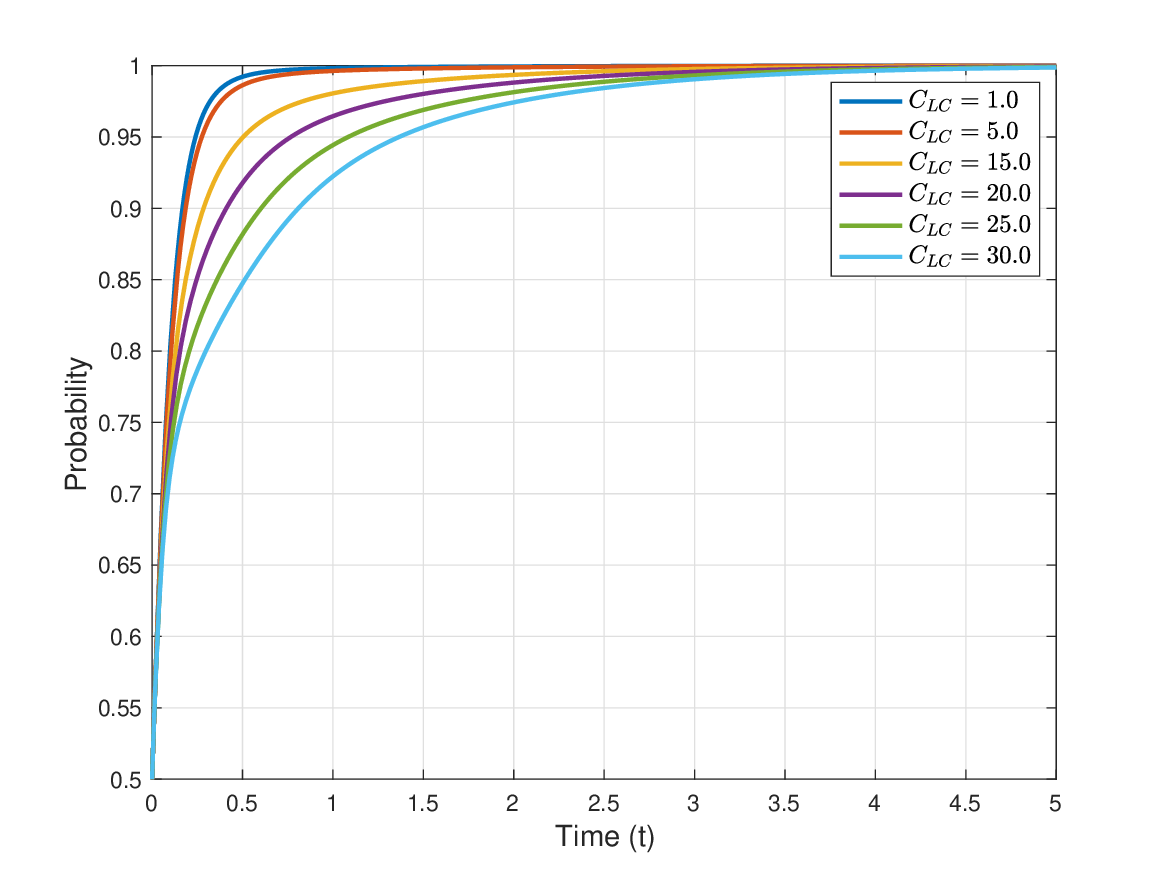}
		\caption{Evolution of the Japanese Fisheries Association over Time.}
		\label{f3c}
	\end{subfigure}
	\caption{Impact of Litigation Compensation to Other Countries $C_{LC}$ on Evolution.}
	\label{f3}
\end{figure}

We analyze the impact of litigation compensation to other countries, denoted as \(C_{LC}\), on the evolutionary outcomes and trajectories of the three game entities. We assume \(C_{LC} = 1.0, 5.0, 15.0, 20.0, 25.0,\) and \(30.0\), while keeping other parameters constant as shown in Table~\ref{tab:initial-parameters}. The specific impacts of varying \(C_{LC}\) are illustrated in Fig.~\ref{f3}.
Fig.~\ref{f3a} demonstrates that as \(C_{LC}\) increases, the evolutionary stable point of the Japanese government, \(\gamma_4(0, 0, 1)^T\), remains constant, but the evolutionary speed slows down, indicating that the Japanese government takes a longer time to reach the equilibrium point. This could signify that higher litigation costs from other countries impose greater economic pressure and uncertainty on Japan, necessitating a more cautious and deliberate approach in strategy formulation.
As depicted in Fig.~\ref{f3b}, with the increase in \(C_{LC}\), the evolutionary stable point of other countries, \(\gamma_4(0, 0, 1)^T\), also remains constant, and the evolutionary speed gradually decelerates. This could be interpreted as other countries facing an increased number of legal disputes and litigations, requiring more resources and time to address these legal matters, thereby inducing greater economic pressure and uncertainty on their own economies.
Fig.~\ref{f3c} illustrates that with the increase in \(C_{LC}\), the evolutionary stable point of the Japanese Fisheries Association, \(\gamma_4(0, 0, 1)^T\), also remains constant, but the evolutionary speed progressively slows down, meaning the Japanese Fisheries Association reaches the equilibrium point more slowly. This could be attributed to the Japanese Fisheries Association needing more time to assess and adapt to the potential impacts of higher litigation compensations paid by the Japanese government to other countries on the fisheries and related industries in Japan.

In conclusion, the varying levels of litigation compensation to other countries not only influence the strategic deliberation and evolutionary pace of the Japanese government and other countries but also reflect the adaptive responses and strategic considerations of the Japanese Fisheries Association in the face of increased economic pressures and uncertainties.

\subsubsection*{Litigation compensation to fisheries association \(C_{LF}\)}

\begin{figure}[H]
	\centering
	\begin{subfigure}[t]{0.32\linewidth}
		\includegraphics[width=\linewidth]{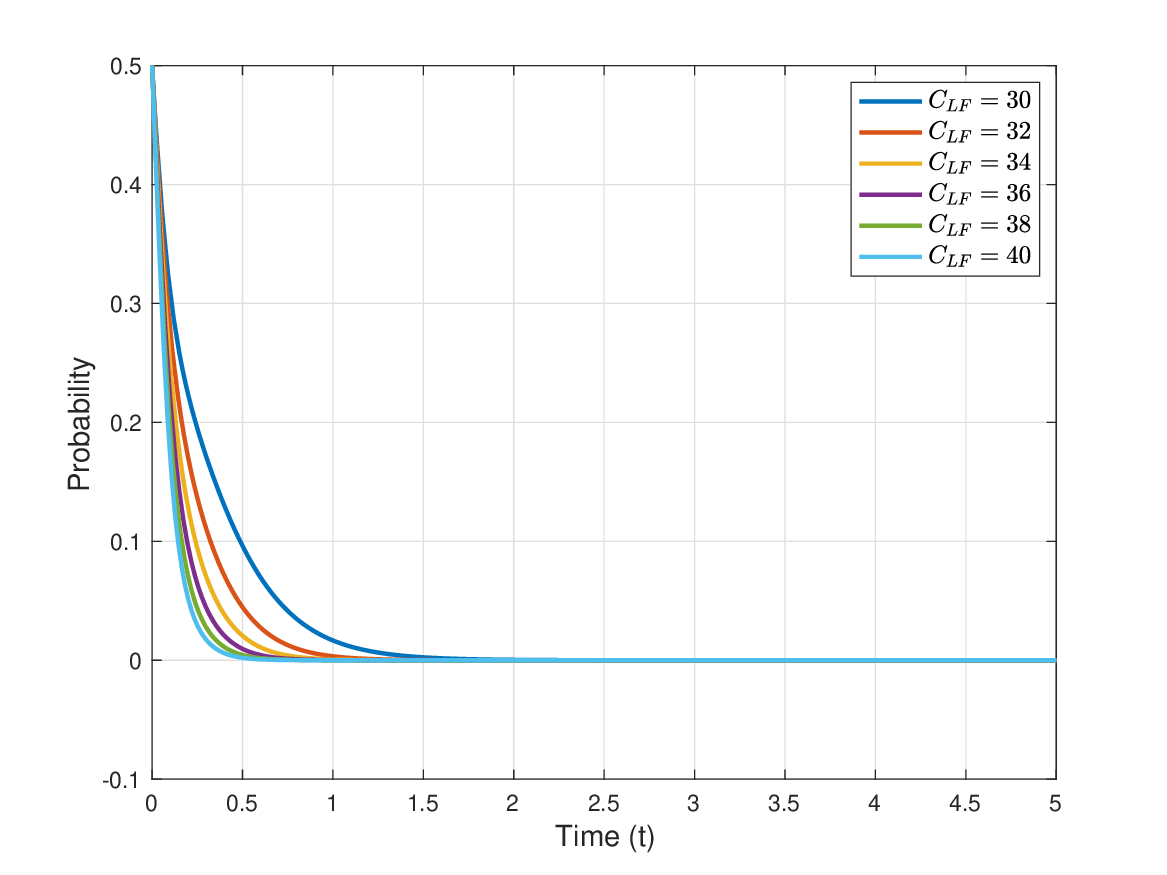}
		\caption{Evolution of the Japanese Government over Time.}
		\label{f4a}
	\end{subfigure}
	\begin{subfigure}[t]{0.32\linewidth}
		\includegraphics[width=\linewidth]{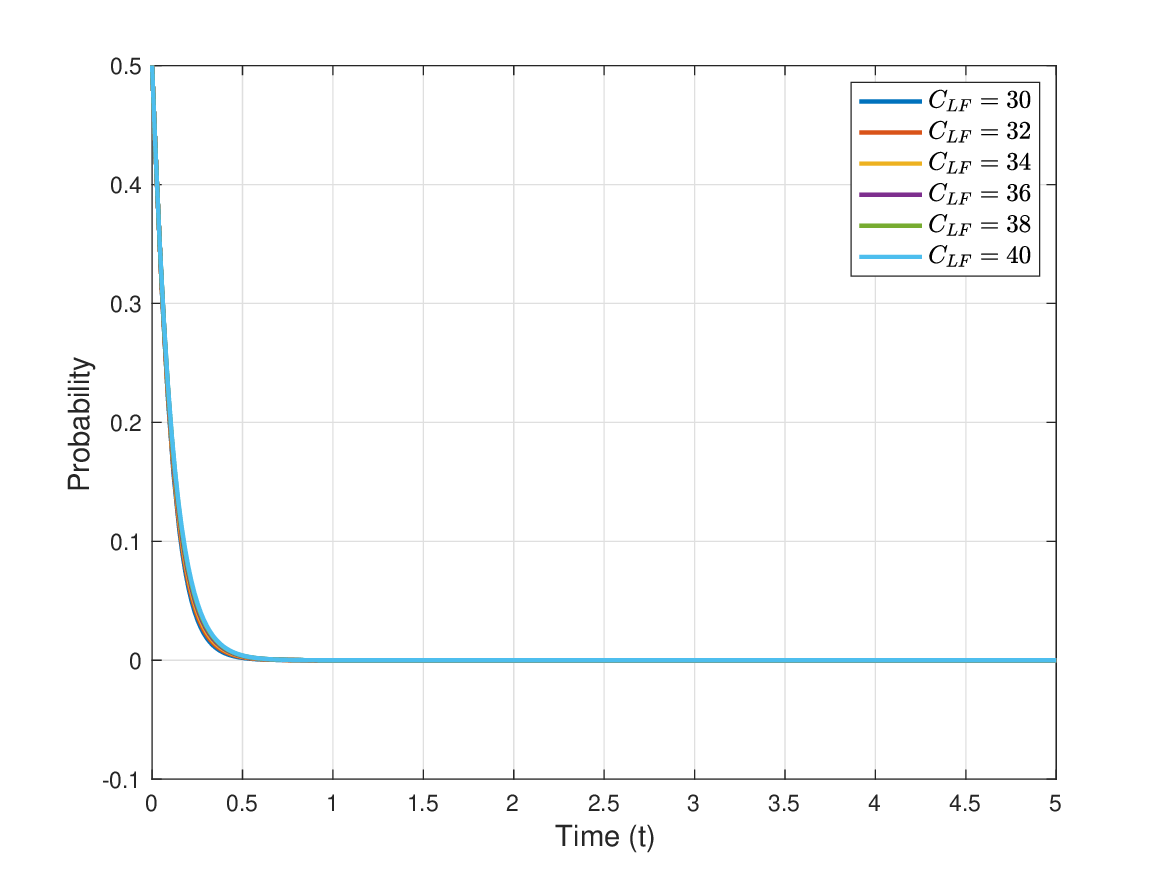}
		\caption{Evolution of Other Countries over Time.}
		\label{f4b}
	\end{subfigure}
	\begin{subfigure}[t]{0.32\linewidth}
		\includegraphics[width=\linewidth]{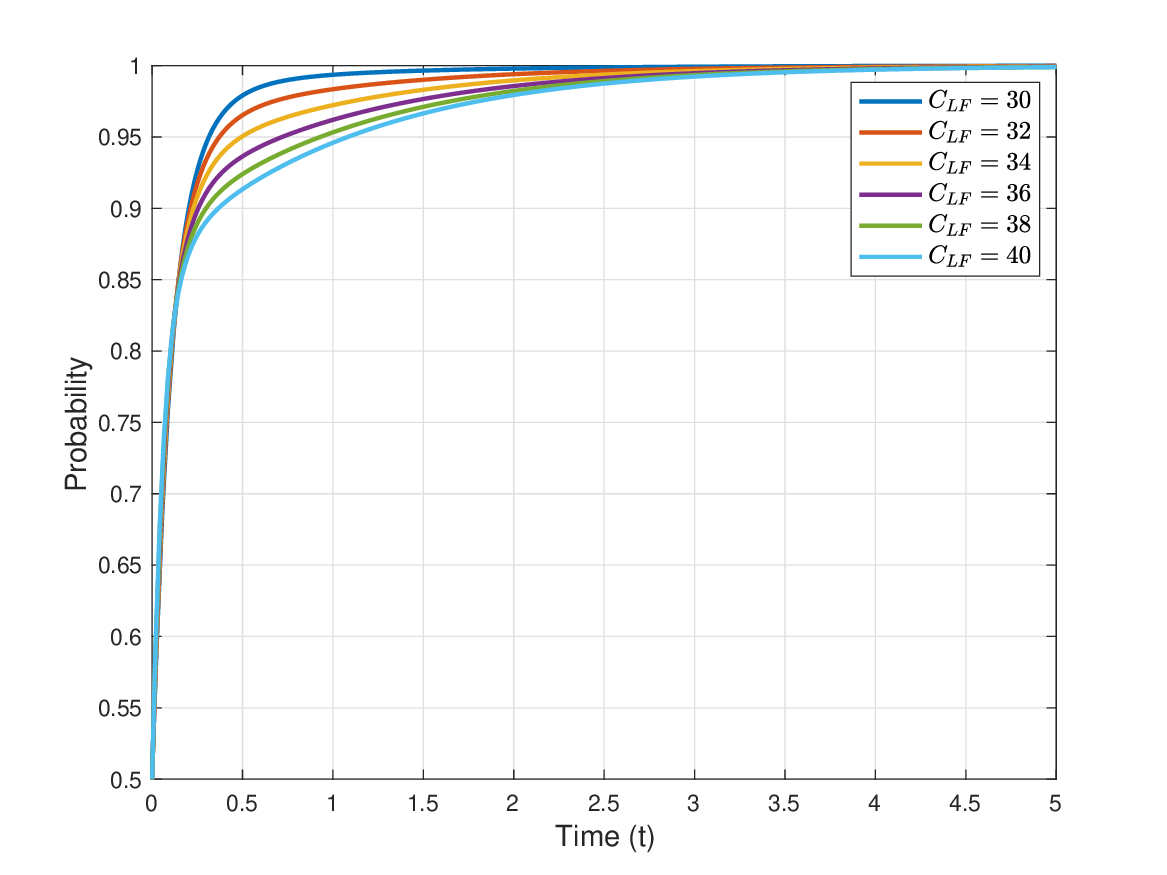}
		\caption{Evolution of the Japanese Fisheries Association over Time.}
		\label{f4c}
	\end{subfigure}
	\caption{Impact of Litigation Compensation to Fisheries Association $C_{LF}$ on Evolution}
	\label{f4}
\end{figure}

We analyze the impact of litigation compensation to the Fisheries Association, denoted as \(C_{LF}\), on the evolutionary outcomes and trajectories of the three game entities. We assume \(C_{LF} = 30, 32, 34, 36, 38,\) and \(40\), while keeping other parameters constant as shown in Table~\ref{tab:initial-parameters}. The specific impacts of varying \(C_{LF}\) are illustrated in Fig.~\ref{f4}.
Fig.~\ref{f4a} demonstrates that as \(C_{LF}\) increases, the evolutionary stable point of the Japanese government, \(\gamma_4(0, 0, 1)^T\), remains constant, but the evolutionary speed rapidly accelerates, meaning the Japanese government reaches the equilibrium point more swiftly. This could signify that higher litigation costs to the Fisheries Association prompt the Japanese government to quickly solidify its strategy to mitigate additional costs induced by uncertainties, emphasizing the government's strategic adaptability in response to internal compensatory obligations.
As depicted in Fig.~\ref{f4b}, with the increase in \(C_{LF}\), the evolutionary stable point of other countries, \(\gamma_4(0, 0, 1)^T\), also remains constant, and the evolutionary speed does not significantly change. This could be interpreted as other countries' decision-making being predominantly influenced by their domestic factors and international political elements, having little association with the compensations paid by the Japanese government to its domestic Fisheries Association, reflecting the relative independence of international actors from Japan's internal compensatory dynamics.
Fig.~\ref{f4c} illustrates that with the increase in \(C_{LF}\), the evolutionary stable point of the Japanese Fisheries Association, \(\gamma_4(0, 0, 1)^T\), also remains constant, but the evolutionary speed progressively slows down, meaning the Japanese Fisheries Association reaches the equilibrium point more slowly. This could be attributed to the Japanese Fisheries Association needing more time to assess and adapt to the potential impacts of receiving higher litigation compensations. Despite the potential influx of compensatory funds, the association must weigh the long-term considerations for the fishermen and the future development of fisheries livelihoods, highlighting the intricate balance between immediate financial relief and long-term sustainable development.

In conclusion, the varying levels of litigation compensation to the Fisheries Association not only influence the strategic deliberation and evolutionary pace of the Japanese government and the Japanese Fisheries Association but also underscore the nuanced interplay between economic compensations and sustainable strategic considerations in the context of environmental disputes.

\subsubsection*{Aid received by Japan from other countries \(C_{HJ}\)}

\begin{figure}[H]
	\centering
	\begin{subfigure}[t]{0.32\linewidth}
		\includegraphics[width=\linewidth]{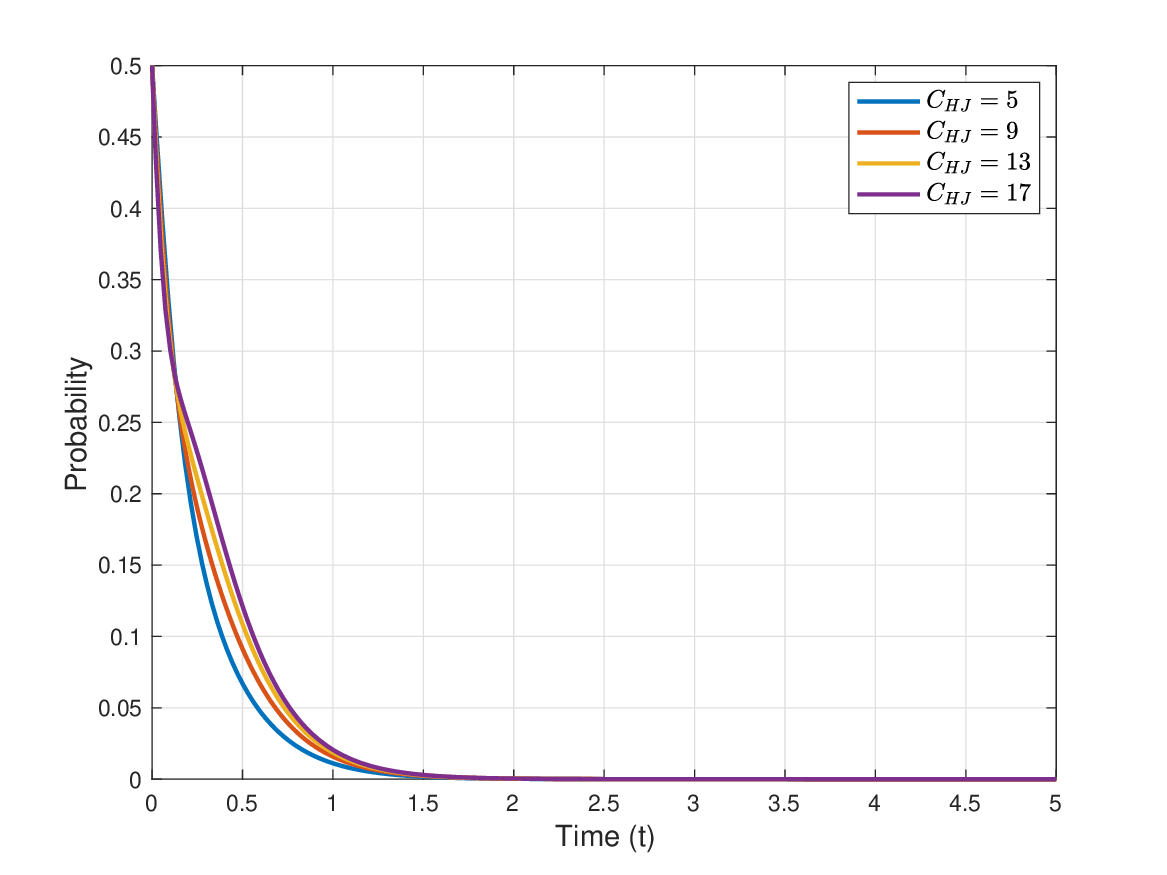}
		\caption{Evolution of the Japanese Government over Time.}
		\label{f5a}
	\end{subfigure}
	\begin{subfigure}[t]{0.32\linewidth}
		\includegraphics[width=\linewidth]{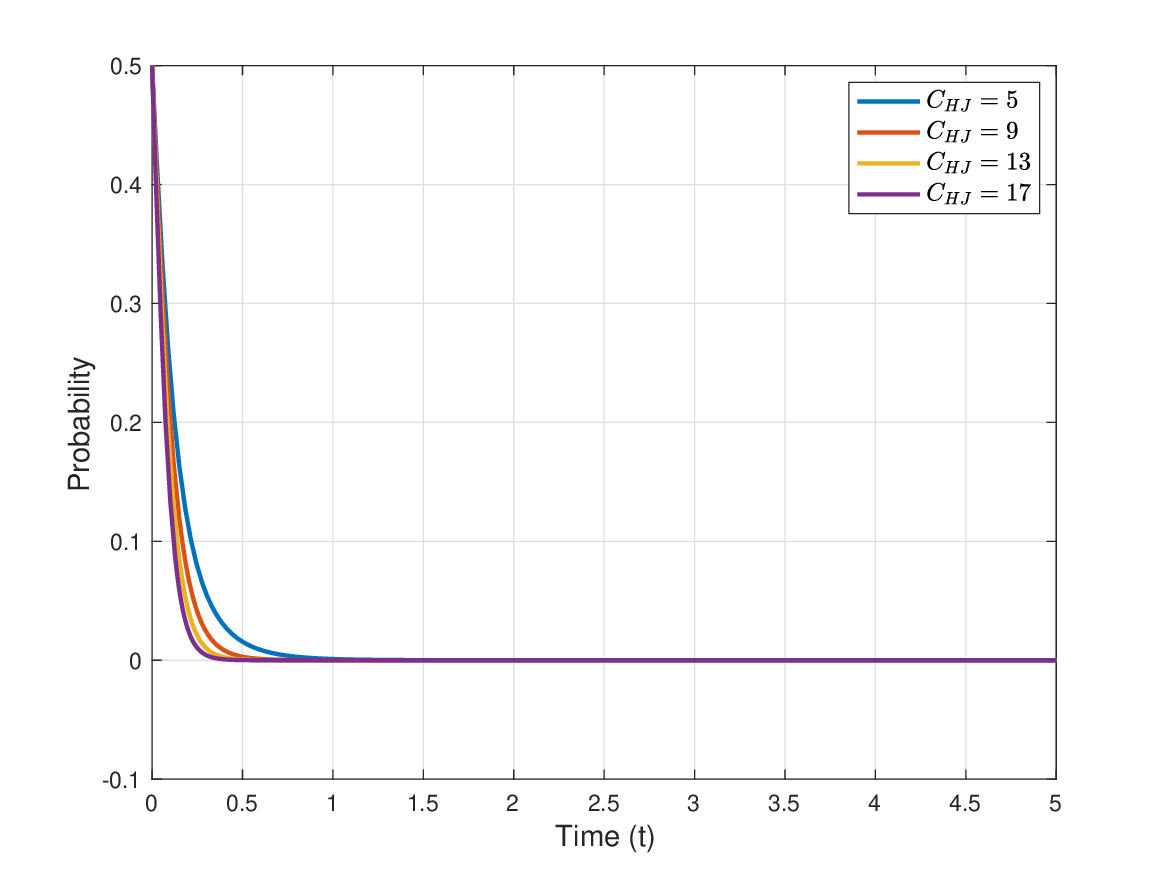}
		\caption{Evolution of Other Countries over Time}
		\label{f5b}
	\end{subfigure}
	\begin{subfigure}[t]{0.32\linewidth}
		\includegraphics[width=\linewidth]{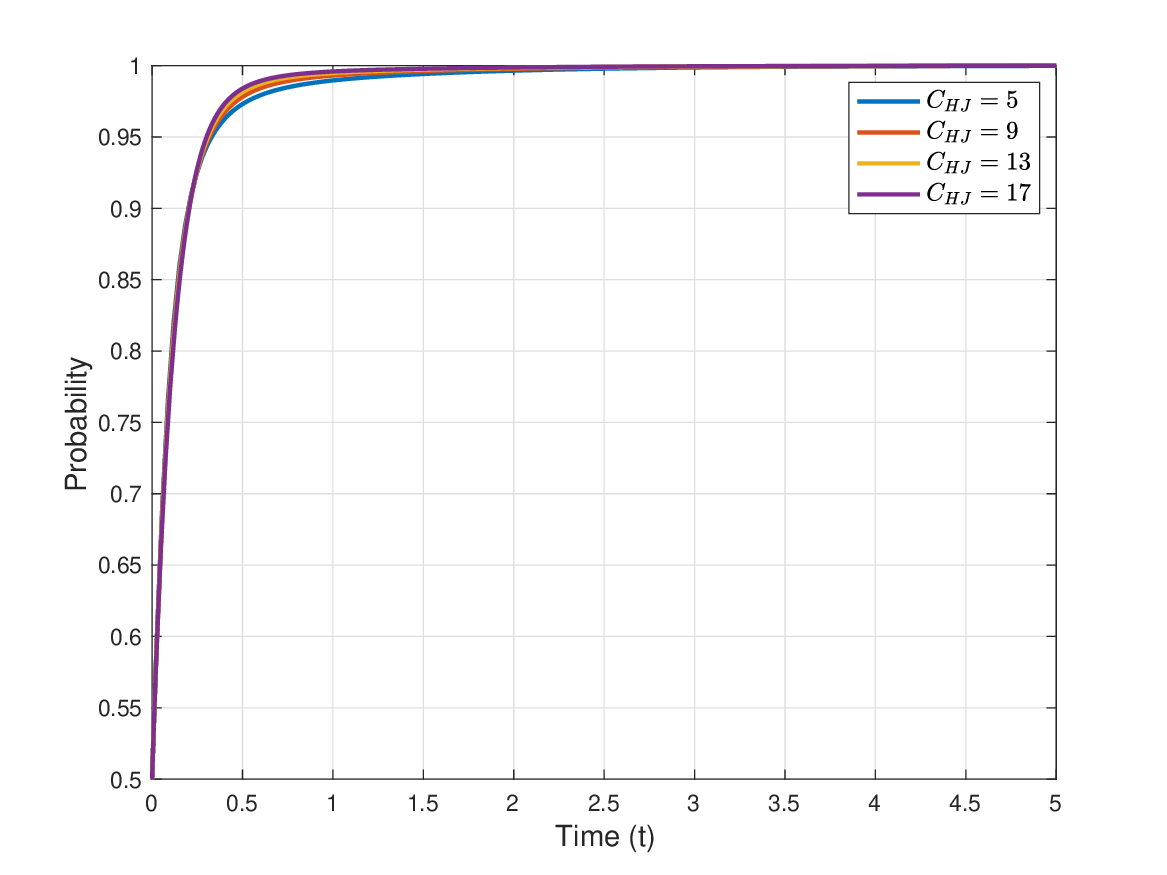}
		\caption{Evolution of the Japanese Fisheries Association over Time.}
		\label{f5c}
	\end{subfigure}
	\caption{Impact of Aid Received by Japan from Other Countries $C_{HJ}$ on Evolution.}
	\label{f5}
\end{figure}

We explore the implications of the aid received by Japan from other countries, denoted as \(C_{HJ}\), on the evolutionary outcomes and trajectories of the three game entities. We assume \(C_{HJ} = 5, 9, 13,\) and \(17\), while maintaining other parameters as shown in Table~\ref{tab:initial-parameters}. The specific impacts of varying \(C_{HJ}\) are illustrated in Fig.~\ref{f5}.
Fig.~\ref{f5a} demonstrates that as \(C_{HJ}\) increases, the evolutionary stable point of the Japanese government, \(\gamma_4(0, 0, 1)^T\), remains constant, but the evolutionary speed progressively decelerates, meaning the Japanese government reaches the equilibrium point more slowly. This suggests that the aid from other countries may induce the Japanese government to deliberate more meticulously on its strategies, with the slower evolutionary speed aimed at mitigating additional costs induced by uncertainties. Enhanced international aid from other countries acts as a deterrent to Japan’s discharge strategy and prolongs the time to reach the evolutionary stable point, serving as an incentivizing measure for Japan to refrain from discharging.
As depicted in Fig.~\ref{f5b}, with the increase in \(C_{HJ}\), the evolutionary stable point of other countries, \(\gamma_4(0, 0, 1)^T\), also remains constant, but the evolutionary speed accelerates. This could be interpreted as the larger the international aid provided by other countries, the slower the progression of the crisis, and when improvements are observed, other countries also accelerate their evolution, reflecting the adaptive nature of international actors in response to the evolving crisis dynamics.
Fig.~\ref{f5c} illustrates that with the increase in \(C_{HJ}\), the evolutionary stable point of the Japanese Fisheries Association, \(\gamma_4(0, 0, 1)^T\), also remains constant, but the evolutionary speed accelerates, meaning the Japanese Fisheries Association reaches the equilibrium point more swiftly. This could be attributed to the Japanese Fisheries Association perceiving the aid from other countries as an incentive that could potentially foster the development of Japanese fisheries and secure the income of fishermen, highlighting the perceived positive impact of international aid on domestic industries.

In conclusion, the varying levels of international aid received by Japan not only influence the strategic deliberation and evolutionary pace of the involved entities but also underscore the intricate interplay between international cooperation and strategic environmental decision-making in the context of international environmental disputes.

\subsubsection*{Initial strategy selection probabilities of Japan, other countries, and the japanese fisheries Association: \(x_0\), \(y_0\), \(z_0\)}

\begin{figure}[H]
	\centering
	\begin{subfigure}[t]{0.32\linewidth}
		\includegraphics[width=\linewidth]{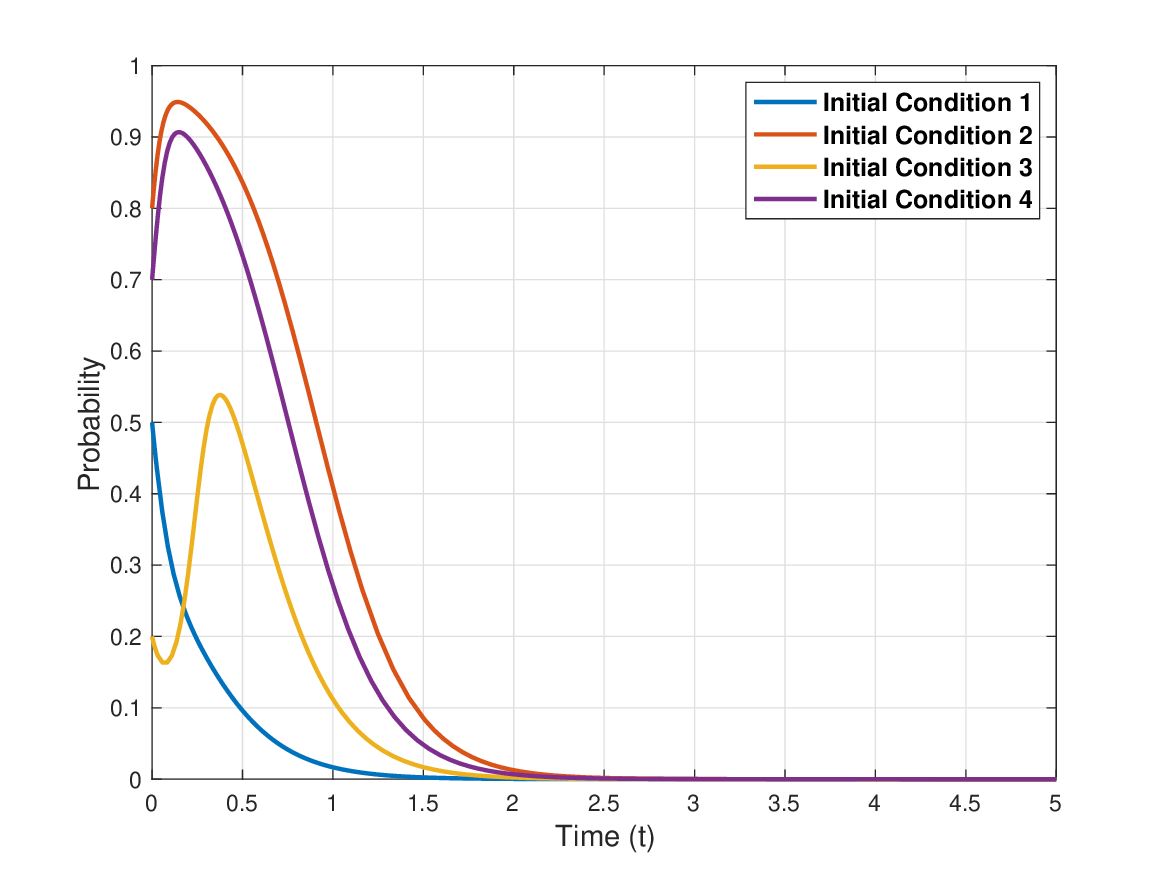}
		\caption{Evolution of the Japanese Government over Time.}
		\label{f6a}
	\end{subfigure}
	\begin{subfigure}[t]{0.32\linewidth}
		\includegraphics[width=\linewidth]{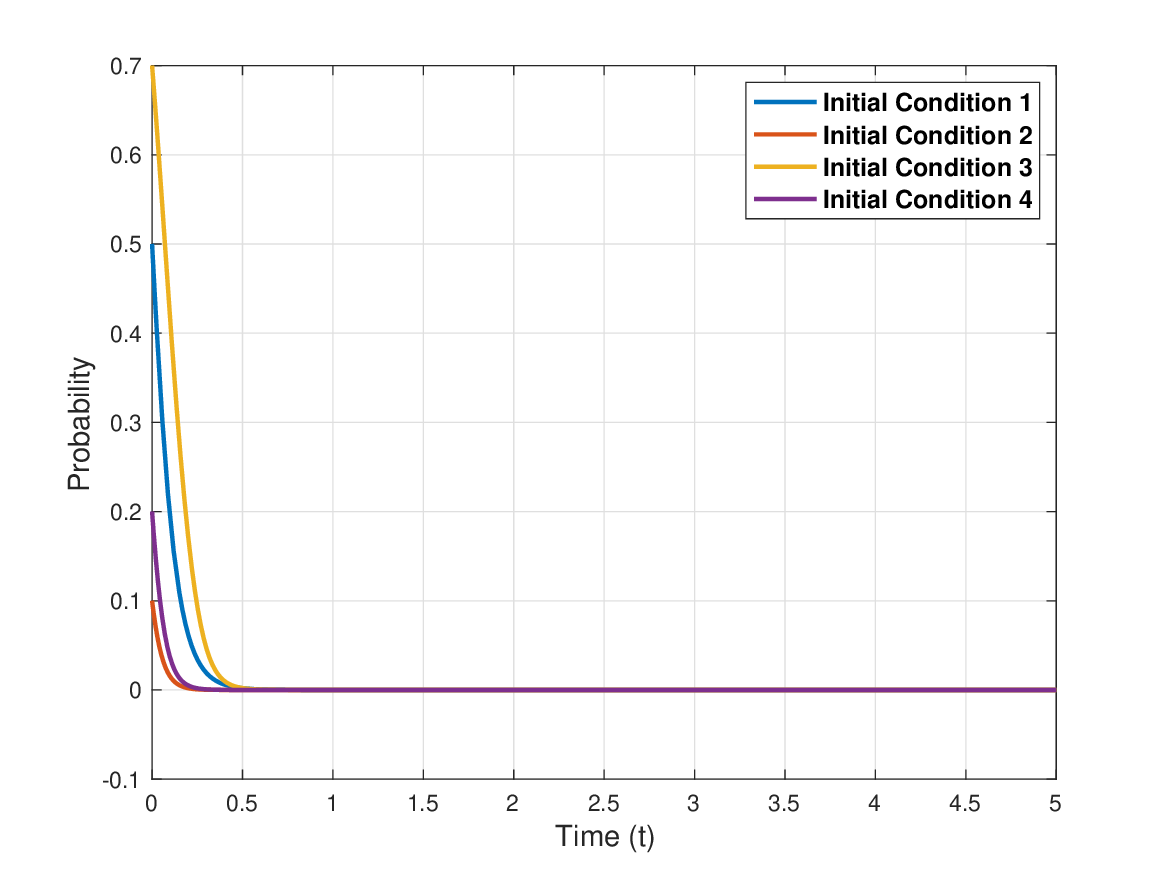}
		\caption{Evolution of Other Countries over Time.}
		\label{f6b}
	\end{subfigure}
	\begin{subfigure}[t]{0.32\linewidth}
		\includegraphics[width=\linewidth]{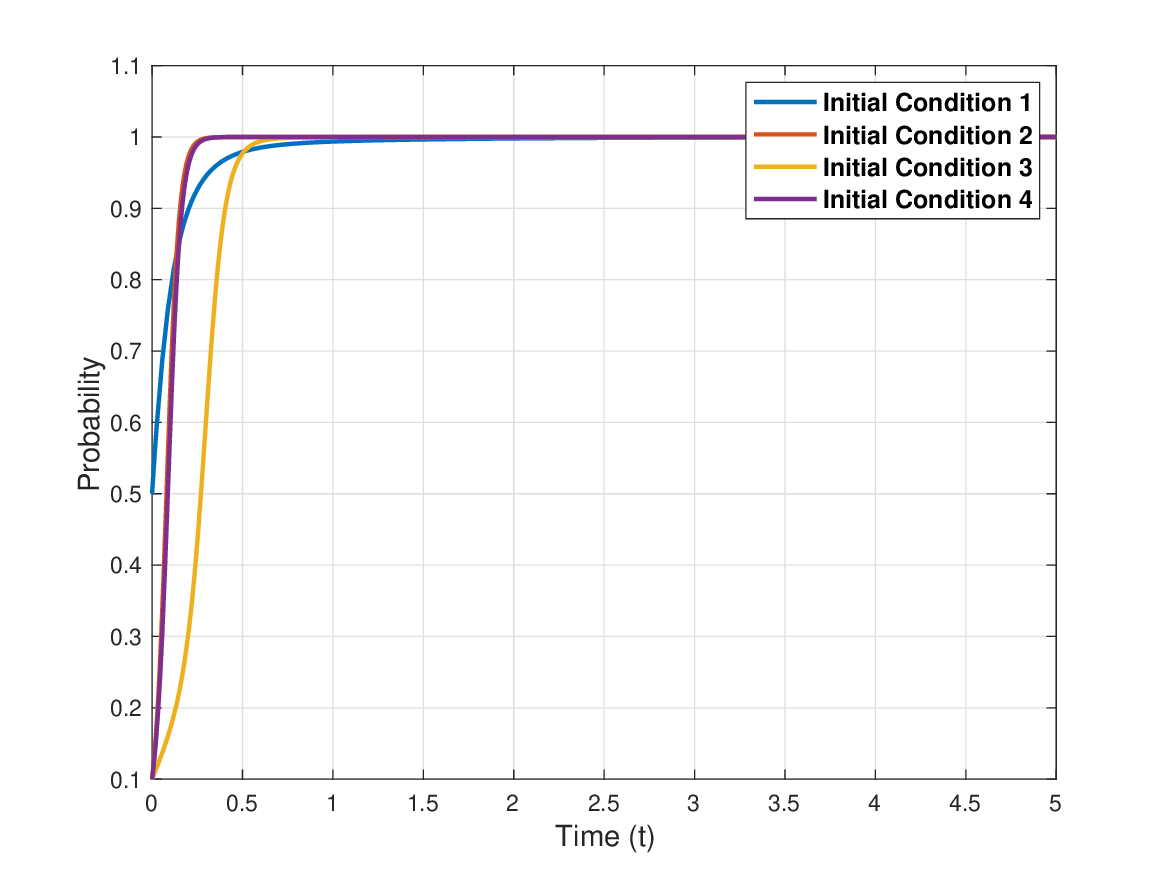}
		\caption{Evolution of the Japanese Fisheries Association over Time.}
		\label{f6c}
	\end{subfigure}
	\caption{Impact of the Initial Strategy Selection Probabilities \((x_0\), \(y_0\), \(z_0)\) of the Three Parties on Evolution.}
	\label{f6}
\end{figure}

Assuming \((x_0, y_0, z_0) = (0.5, 0.5, 0.5), (0.8, 0.1, 0.1), (0.2, 0.7, 0.1), (0.7, 0.2, 0.1)\), while maintaining other parameters as shown in Table~\ref{tab:initial-parameters}, we can observe the evolutionary outcomes and trajectories of the three game entities under varying initial strategy selection probabilities, as depicted in Fig.~\ref{f6}. 
From the figure, it is evident that when initial strategy selection probabilities vary, the evolutionary stable point \(\gamma_4(0, 0, 1)^T\) remains constant, indicating that it does not affect the evolutionary outcomes but does influence the evolutionary trajectories. For Japan, a stronger initial strategy selection probability, \(x_0\), will delay Japan's discharge strategy and decelerate its evolutionary speed, as illustrated in Fig.~\ref{f6a}. 
For other countries, the evolutionary stable point \(\gamma_4(0, 0, 1)^T\) also remains constant when initial strategy selection probabilities vary, signifying that it does not impact the evolutionary outcomes but does have certain effects on the evolutionary trajectories. A stronger initial strategy selection probability, \(y_0\), will delay the opposition strategy of the Fisheries Association and slow down its evolutionary speed, as shown in Fig.~\ref{f6b}.
Lastly, for the Fisheries Association, the evolutionary stable point \(\gamma_4(0, 0, 1)^T\) still remains constant when initial strategy selection probabilities change, suggesting that it does not alter the evolutionary outcomes but does affect the evolutionary trajectories. A stronger initial strategy selection probability, \(z_0\), will delay the opposition strategy of the Fisheries Association and decelerate its evolutionary speed, as depicted in Fig.~\ref{f6c}.

\section*{Conclusions and policy implications}
\subsection*{Conclusions}

To accurately represent the real-world issue of Japan's nuclear wastewater discharge within the model and to explore the Evolutionarily Stable Strategies (ESS) of Japan, other countries, and the Japanese Fisheries Association, this study employs an evolutionary game model. This model analyzes the strategic interactions and evolutionary dynamics of the three entities on the issue of nuclear wastewater discharge, identifying three ESSs and their respective stability conditions. Before conducting numerical simulations, assumptions and stability analyses were made on the parameters. We set a group of values based on stability conditions and obtained the evolutionary process through numerical simulations. Finally, we conducted an in-depth analysis of the impact generated by some key parameters on the stable point where Japan ceases discharge. The innovation of this study primarily lies in exploring the stable situation where Japan does not discharge and studying the impact of different parameters on evolutionary trajectories. Especially against the backdrop of Japan commencing nuclear wastewater discharge on August 24, 2023, our research is novel and practically significant, providing important theoretical support for sustainability and environmental protection. The main findings of this study are as follows:

Firstly, we distinguished between the three ESSs and their respective stability conditions, identifying three ESSs and their corresponding stability conditions. Secondly, we determined the impact of key parameters on evolutionary outcomes and trajectories. Japan's international image, litigation compensation from other countries, reduction in export tax revenue for Japan due to discharge, international aid, litigation compensation to the Fisheries Association, discharge costs, marine monitoring costs, nuclear wastewater storage costs, the international image of the Fisheries Association, potential gains of other countries, and additional costs of other countries are key factors affecting the evolutionary outcomes of the three entities.

We found that variations in these parameters would influence the evolutionary speed and strategy choices of each entity, but due to the assumptions made earlier, the evolutionary stable points remained constant. The strategic choices of Japan depend on the comparison between discharge costs and litigation costs. For instance, higher discharge costs and litigation compensation would prompt the Japanese government to solidify its strategy more swiftly to reduce additional costs brought by uncertainties. Conversely, international aid might induce the Japanese government to deliberate more on its strategies, slowing down the evolutionary speed. When Japan chooses the discharge strategy, the choices of other countries are determined by the trade-off between litigation compensation and costs, while the choices of the Fisheries Association are based on the comparison between litigation compensation, litigation costs, and the image of the Association.

\subsection*{Policy implications}

This paper proposes the following policy implications from the perspectives of Japan, other countries, and the Japanese Fisheries Association to promote Japan's strategic transition from discharging to not discharging. Additionally, it analyzes the maximum benefits that other countries and the Fisheries Association can obtain in the scenario where Japan persists in its discharge strategy.

Additionally, it may be prudent for Japan to reevaluate the viability of existing alternative treatment solutions and deliberate over the discontinuation of sea discharges. The considerations should extend beyond mere cost-benefit analysis; the profound detriment to marine ecosystems, substantial litigation compensations, prospective forfeitures of international aid, and acute interest conflicts with other nations and the Fisheries Association necessitate meticulous attention from Japan. Presently, a spectrum of technologies and apparatus, including catalytic exchange, cryogenic distillation, and water distillation technologies, are available, capable of effectively purifying various radioactive nuclides in nuclear-contaminated water, showcasing the potential to eliminate tritium and other nuclides. Consequently, Japan should refrain from selecting the dilution discharge option merely due to its cost-effectiveness. It is unequivocally irresponsible to discharge nuclear wastewater into the sea, as it imposes enduring risks to the global ecosystem and humanity at large.
On the other hand, a consensus should be reached among other countries and the Japanese Fisheries Association to staunchly oppose and sanction Japan's discharge strategy, especially when the costs of compensation are comparatively low. Neighboring countries of Japan should not be burdened with the enduring and substantial costs of marine monitoring and cleansing resulting from Japan's discharge. It is imperative for these nations to promptly safeguard their marine and economic environments against the repercussions of such discharge.

The Fisheries Association, as the primary victim of Japan's discharge strategy, representing the interests of the vast number of Japanese fishermen, should more firmly adopt an opposing stance, staunchly defending life, health, and livelihoods. For instance, the Japanese Fisheries Association can intensify the opposition to the nuclear wastewater discharge plan and strive for more support and attention to protect the sustainable development of fisheries and the ecological environment.

\ \\
For countries neighboring Japan:
\begin{itemize}
	\item Increase Diplomatic Pressure: Express concerns and opposition to nuclear wastewater discharge to the Japanese government through diplomatic channels. Collaborate with other neighboring countries to form a united stance and jointly exert diplomatic pressure on Japan.
	\item Collaborate on Scientific Research: Cooperate with scientists, research institutions, and fisheries departments of neighboring countries to jointly conduct scientific research to assess the potential impacts of nuclear wastewater discharge on surrounding marine areas and fisheries resources. Strengthen the scientific basis of opposition plans through scientific data and expert opinions.
	\item Enhance Propaganda and Education: Publicize the risks and impacts of nuclear wastewater discharge through media and educational channels. Intensify propaganda, raise public environmental awareness and attention to the issue of nuclear wastewater discharge, and create public opinion pressure.
\end{itemize}

\noindent For countries far from Japan:
\begin{itemize}
	\item Issue Statements and Express Concerns: Express concerns and opposition to nuclear wastewater discharge to the Japanese government through diplomatic channels. International organizations, the United Nations, and other international platforms can also be utilized to express concerns and opposition.
	\item Strengthen International Cooperation: Monitor and assess Japan's nuclear wastewater discharge plan jointly with other countries through international cooperation and organizations. Form a stronger international opposition voice through sharing scientific data and expert opinions.
	\item Provide Technical Support and Expert Opinions: Offer technical support and expert opinions to countries neighboring Japan to help them assess the potential impacts of nuclear wastewater discharge on their fisheries and environment. Strengthen cooperation and exchange to form a closer international cooperation network. These measures can help other countries strengthen the opposition to Japan's nuclear wastewater discharge plan and jointly maintain the sustainable development of the marine environment and fisheries resources.
\end{itemize}

\section*{Appendix}
\subsection*{Jacobian matrices for different points}
\begin{equation}
	\boldsymbol{J}(0,0,0)=\left[\begin{array}{ccc}
		\mathrm{C}_{\mathrm{SJ}}-\mathrm{C}_{\mathrm{MJ}}-\mathrm{C}_{\mathrm{DJ}}, & 0, & 0 \\
		0, & -\mathrm{C}_{\mathrm{HJ}}, & 0 \\
		0, & 0, & \mathrm{C}_{\mathrm{IF}}
	\end{array}\right]
\end{equation}

\begin{equation}
	\boldsymbol{J}(0,0,1)=\left[\begin{array}{ccc}
		\mathrm{C}_{\mathrm{SJ}}-\mathrm{C}_{\mathrm{LF}}-\mathrm{C}_{\mathrm{MJ}}-\mathrm{C}_{\mathrm{DJ}}, & 0, & 0 \\
		0, & -\mathrm{C}_{\mathrm{HJ}}, & 0 \\
		0, & 0, & -\mathrm{C}_{\mathrm{IF}}
	\end{array}\right]
\end{equation}

\begin{equation}
	J(0,1,0)=\left[\begin{array}{ccc}
		\mathrm{C}_{\mathrm{SJ}}-\mathrm{C}_{\mathrm{HJ}}-\mathrm{C}_{\mathrm{LC}}-\mathrm{C}_{\mathrm{MJ}}-\mathrm{C}_{\mathrm{DJ}}-\mathrm{I}_{\mathrm{J}}-\mathrm{T}_{\mathrm{RJ}}, & 0, & 0 \\
		0, & \mathrm{C}_{\mathrm{HJ}}, & 0 \\
		0, & 0, & \mathrm{C}_{\mathrm{IF}}
	\end{array}\right]
\end{equation}

\begin{equation}
	\boldsymbol{J}(0,1,1)=\left[\begin{array}{ccc}
		\mathrm{C}_{\mathrm{SJ}}-\mathrm{C}_{\mathrm{HJ}}-\mathrm{C}_{\mathrm{LC}}-\mathrm{C}_{\mathrm{LF}}-\mathrm{C}_{\mathrm{MJ}}-\mathrm{C}_{\mathrm{DJ}}-\mathrm{I}_{\mathrm{J}}-\mathrm{T}_{\mathrm{RJ}}, & 0, & 0 \\
		0, & \mathrm{C}_{\mathrm{HJ}}, & 0 \\
		0, & 0, & -\mathrm{C}_{\mathrm{IF}}
	\end{array}\right]
\end{equation}

\begin{equation}
	\boldsymbol{J}(1,0,0)=\left[\begin{array}{ccc}
		\mathrm{C}_{\mathrm{DJ}}+\mathrm{C}_{\mathrm{MJ}}-\mathrm{C}_{\mathrm{SJ}}, & 0, & 0 \\
		0, & \mathrm{~B}_{\mathrm{SP}}+\mathrm{C}_{\mathrm{LC}}-\mathrm{C}_{\mathrm{SC}}, & 0 \\
		0, & 0, & \mathrm{C}_{\mathrm{LF}}+\mathrm{C}_{\mathrm{IF}}
	\end{array}\right]
\end{equation}

\begin{equation}
	\boldsymbol{J}(1,0,1)=\left[\begin{array}{ccc}
		\mathrm{C}_{\mathrm{DJ}}+\mathrm{C}_{\mathrm{LF}}+\mathrm{C}_{\mathrm{MJ}}-\mathrm{C}_{\mathrm{SJ}}, & 0, & 0 \\
		0, & \mathrm{~B}_{\mathrm{SP}}+\mathrm{C}_{\mathrm{LC}}-\mathrm{C}_{\mathrm{SC}}, & 0 \\
		0, & 0, & -\mathrm{C}_{\mathrm{IF}}-\mathrm{C}_{\mathrm{LF}}
	\end{array}\right]
\end{equation}

\begin{equation}
	\boldsymbol{J}(1,1,0)=\left[\begin{array}{ccc}
		\mathrm{C}_{\mathrm{DJ}}+\mathrm{C}_{\mathrm{HJ}}+\mathrm{C}_{\mathrm{LC}}+\mathrm{C}_{\mathrm{MJ}}-\mathrm{C}_{\mathrm{SJ}}+\mathrm{I}_{\mathrm{J}}+\mathrm{T}_{\mathrm{RJ}}, & 0, & 0 \\
		0, & \mathrm{C}_{\mathrm{SC}}-\mathrm{C}_{\mathrm{LC}}-\mathrm{B}_{\mathrm{SP}} & 0 \\
		0, & 0, & \mathrm{C}_{\mathrm{LF}}+\mathrm{C}_{\mathrm{IF}}
	\end{array}\right]
\end{equation}

\begin{equation}
	J(1,1,1)=\left[\begin{array}{ccc}
		\mathrm{C_{DJ}}+\mathrm{C_{HJ}}+\mathrm{C_{LC}}+\mathrm{C_{LF}}+\mathrm{C_{MJ}}-\mathrm{C_{SJ}}+\mathrm{I_J}+\mathrm{T_{RJ}},&0, & 0 \\
		0, &\mathrm{C_{SC}}-\mathrm{C_{LC}}-\mathrm{B_{SP}}, & 0 \\
		0, & 0,&-\mathrm{C_{IF}}-\mathrm{C_{LF}}
	\end{array}\right]
\end{equation}

\subsection*{Eigenvalues for different points}

\begin{equation}
	\gamma_{1}(0,0,0)=\left(\begin{array}{c}
		\mathrm{C}_{\mathrm{SJ}}-\mathrm{C}_{\mathrm{MJ}}-\mathrm{C}_{\mathrm{DJ}} \\
		-\mathrm{C}_{\mathrm{HJ}} \\
		\mathrm{C}_{\mathrm{IF}}
	\end{array}\right)
\end{equation}

\begin{equation}
	\gamma_{2}(1,0,0)=\left(\begin{array}{c}
		\mathrm{C}_{\mathrm{DJ}}+\mathrm{C}_{\mathrm{MJ}}-\mathrm{C}_{\mathrm{SJ}} \\
		\mathrm{B}_{\mathrm{SP}}+\mathrm{C}_{\mathrm{LC}}-\mathrm{C}_{\mathrm{SC}} \\
		\mathrm{C}_{\mathrm{LF}}+\mathrm{C}_{\mathrm{IF}}
	\end{array}\right)
\end{equation}

\begin{equation}
	\gamma_{3}(0,1,0)=\left(\begin{array}{c}
		\mathrm{C}_{\mathrm{SJ}}-\mathrm{C}_{\mathrm{HJ}}-\mathrm{C}_{\mathrm{LC}}-\mathrm{C}_{\mathrm{MJ}}-\mathrm{C}_{\mathrm{DJ}}-\mathrm{I}_{\mathrm{J}}-\mathrm{T}_{\mathrm{RJ}} \\
		\mathrm{C}_{\mathrm{HJ}} \\
		\mathrm{C}_{\mathrm{IF}}
	\end{array}\right)
\end{equation}

\begin{equation}
	\gamma_{4}(0,0,1)=\left(\begin{array}{c}
		\mathrm{C}_{\mathrm{SJ}}-\mathrm{C}_{\mathrm{LF}}-\mathrm{C}_{\mathrm{MJ}}-\mathrm{C}_{\mathrm{DJ}} \\
		-\mathrm{C}_{\mathrm{HJ}} \\
		-\mathrm{C}_{\mathrm{IF}}
	\end{array}\right)
\end{equation}

\begin{equation}
	\gamma_{5}(1,1,0)=\left(\begin{array}{c}
		\mathrm{C}_{\mathrm{DJ}}+\mathrm{C}_{\mathrm{HJ}}+\mathrm{C}_{\mathrm{LC}}+\mathrm{C}_{\mathrm{MJ}}-\mathrm{C}_{\mathrm{SJ}}+\mathrm{I}_{\mathrm{J}}+\mathrm{T}_{\mathrm{RJ}} \\
		\mathrm{C}_{\mathrm{SC}}-\mathrm{C}_{\mathrm{LC}}-\mathrm{B}_{\mathrm{SP}} \\
		\mathrm{C}_{\mathrm{LF}}+\mathrm{C}_{\mathrm{IF}}
	\end{array}\right)
\end{equation}

\begin{equation}
	\gamma_{6}(1,0,1)=\left(\begin{array}{c}
		\mathrm{C}_{\mathrm{DJ}}+\mathrm{C}_{\mathrm{LF}}+\mathrm{C}_{\mathrm{MJ}}-\mathrm{C}_{\mathrm{SJ}} \\
		\mathrm{B}_{\mathrm{SP}}+\mathrm{C}_{\mathrm{LC}}-\mathrm{C}_{\mathrm{SC}} \\
		-\mathrm{C}_{\mathrm{IF}}-\mathrm{C}_{\mathrm{LF}}
	\end{array}\right)
\end{equation}

\begin{equation}
	\gamma_{7}(1,0,1)=\left(\begin{array}{c}
		\mathrm{C}_{\mathrm{SJ}}-\mathrm{C}_{\mathrm{HJ}}-\mathrm{C}_{\mathrm{LC}}-\mathrm{C}_{\mathrm{LF}}-\mathrm{C}_{\mathrm{MJ}}-\mathrm{C}_{\mathrm{DJ}}-\mathrm{I}_{\mathrm{J}}-\mathrm{T}_{\mathrm{RJ}} \\
		\mathrm{C}_{\mathrm{HJ}} \\
		-\mathrm{C}_{\mathrm{IF}}
	\end{array}\right)
\end{equation}

\begin{equation}
	\gamma_{8}(1,1,1)=\left(\begin{array}{c}
		\mathrm{C_{DJ}}+\mathrm{C_{HJ}}+\mathrm{C_{LC}}+\mathrm{C_{LF}}+\mathrm{C_{MJ}}-\mathrm{C_{SJ}}+\mathrm{I_J}+\mathrm{T_RJ} \\
		\mathrm{C_{SC}}-\mathrm{C_{LC}}-\mathrm{B_{SP}} \\
		-\mathrm{C_{IF}}-\mathrm{C_{LF}}
	\end{array}\right)
\end{equation}

\section*{Acknowledgements}
This work is supported by the Natural Science Foundation of Xinjiang Uygur Autonomous Region (Grant No. 2021D01A202), the Youth Foundation of China University of Petroleum-Beijing at Karamay (No.XQZX20230034), and the Xinjiang Uygur Autonomous Region Undergraduate Education Reform Project (Grant No. PT2021083).

\section*{Author contributions}
Mingyang Li performed the Methodology, Writing - Original Draft, Visualization and Software; Han Pengsihua performed the Writing - Original Draft; Songqing Zhao performed the Supervision; Zejun Wang performed the Project Administration, Methodology and Writing - Review \& Editing; Limin Yang performed the Funding Acquisition and Writing - Review \& Editing; Weian Liu performed the Writing - Review \& Editing.

\section*{Competing interests}
Te authors declare no competing interests.

\end{document}